\documentclass[12pt]{article}

\usepackage{amsmath,amssymb,amsthm,amsfonts}
\usepackage{geometry}
\usepackage{amscd}
\theoremstyle{plain}
\newtheorem{thm}{Theorem}
\newtheorem{lemma}[thm]{Lemma}
\newtheorem*{claim}{Claim}
\newtheorem{prop}{Proposition}

\theoremstyle{definition}
\newtheorem*{defn}{Definition}
\newtheorem*{notn}{Notation}

\newcommand{\cov}{\operatorname{cov}}
\newcommand{\re}{\operatorname{Re}}
\newcommand{\im}{\operatorname{Im}}

\title{Bourgain's Entropy Estimates Revisited}

\author{John T. Workman}
\date{}

\begin{document}

\maketitle

\begin{abstract}
This serves as a near complete set of notes on Bourgain's well-known
paper \emph{Almost sure convergence and bounded
entropy}~\cite{bourgain}. The two entropy results are treated, as is
one of the applications.  These notes are designed to be independent
of Bourgain's paper and self-contained. There are, at times,
differences between Bourgain's notation and my own.  The same goes
for organization.  However, the proofs herein are essentially his.
\end{abstract}

\vskip0.5in

\section{Preliminaries}

Our setting will be a probability space $(X,\mathcal{F},\mu)$.  We
are interested in certain sequences of operators.  In particular,
given a sequence of operators on $L^2(\mu)$ we want to make some
uniform estimate on their entropy.  What kind of conditions must
this sequence satisfy?  We are led to the following definitions.

\begin{defn} Given a pseudo-metric space $(Y,d)$ (that is, $d$ need
not separate points) and a subset $S \subset Y$, define the
$\delta$-entropy number of $S$ to be the minimal number of (closed)
$\delta$-balls in the $d$ pseudo-metric needed to cover $S$.  We
denote this by $N(S,d,\delta)$. \end{defn}

\begin{notn} Denote by $\mathcal{M}(X)$ the set of measurable
functions $f : X \rightarrow \mathbb{R}$.  By $L^p(\mu)$, it will
always be meant the subset of $\mathcal{M}(X)$ which has finite
$L^p$-norm.  That is, $L^p(\mu)$ consists of real-valued functions
only.
\end{notn}

\begin{defn} Let $T_j : \mathcal{M}(X) \rightarrow \mathcal{M}(X)$,
$j \in \mathbb{N}$, be a sequence of linear operators.  We say
$(T_j)$ is a Bourgain sequence if the following are satisfied:

\begin{enumerate}
\item $T_j : L^1(\mu) \rightarrow L^1(\mu)$ are bounded,

\item $T_j : L^2(\mu) \rightarrow L^2(\mu)$ are isometries,

\item $T_j$ are positive, i.e., if $f \geq 0$ a.s.[$\mu$] then
$T_j(f) \geq 0$ a.s.[$\mu$],

\item $T_j(1) = 1$ (here, 1 refers to the constant function 1),

\item $T_j$ satisfy a mean ergodic condition, i.e.,

\begin{equation*}
\frac{1}{J} \sum_{j=1}^J T_j f \rightarrow \int_X f(x) \, \mu(dx)
\quad \text{in } L^1(\mu) \text{-norm for all } f \in L^1(\mu).
\end{equation*}
\end{enumerate}
\end{defn}

The mean ergodic condition will prove useful as is.  The first four
assumptions lead to the following properties.

\begin{lemma}\label{lemma:T} Let $T : \mathcal{M}(X) \rightarrow
\mathcal{M}(X)$ be a linear operator which satisfies assumptions (1)
- (4) above.  Then, $T : L^\infty(\mu) \rightarrow L^\infty(\mu)$
with $\|T f\|_{L^\infty(\mu)} \leq \|f\|_{L^\infty(\mu)}$. Further,
$T(f^2) = T(f)^2$ a.s.[$\mu$] for all $f \in L^2(\mu)$.
\end{lemma}

\begin{proof}
By assumption 4, $T(c) = c$ for all constant functions $c$.  By
assumption 3, $T(g) \leq T(h)$ a.s.[$\mu$] whenever $g \leq h$
a.s.[$\mu$]. Let $f \in L^\infty(\mu)$.  Then, $f \leq
\|f\|_{L^\infty(\mu)}$ a.s.[$\mu$], so that $T(f) \leq
T(\|f\|_{L^\infty(\mu)}) = \|f\|_{L^\infty(\mu)}$ a.s.[$\mu$].
Similarly, $-T(f) = T(-f) \leq T(\|f\|_{L^\infty(\mu)}) =
\|f\|_{L^\infty(\mu)}$ a.s.[$\mu$].  Thus, $|T(f)| \leq
\|f\|_{L^\infty(\mu)}$ a.s.[$\mu$], giving the first statement.

We now approach the second statement.  First, suppose $A \in
\mathcal{F}$. As $0 \leq \chi_A \leq 1$, we have $0 \leq T(\chi_A)
\leq 1$. By assumption 2, $\mu(A) = \|\chi_A\|_{L^2(\mu)}^2 =
\|T(\chi_A)\|_{L^2(\mu)}^2 = \int_X T(\chi_A)^2 \, d\mu$. On the
other hand, $\mu(A) = 1 - \|\chi_{A^c}\|_{L^2(\mu)}^2 = 1 - \|1 -
\chi_A\|_{L^2(\mu)}^2 = 1 - \|T(1- \chi_A)\|_{L^2(\mu)}^2 = 1 - \|1
- T(\chi_A)\|_{L^2(\mu)}^2 = 1 - \int_X (1 - T(\chi_A))^2 \, d\mu =
\int_X 2T(\chi_A) - T(\chi_A)^2 \, d\mu$.  Setting these two
expressions of $\mu(A)$ equal, we have $\int_X T(\chi_A) \, d\mu =
\int_X T(\chi_A)^2 \, d\mu$.  As $0 \leq T(\chi_A) \leq 1$
a.s.[$\mu$], it must be that $T(\chi_A) = 0, 1$ a.s.[$\mu$]. Namely,
$T$ takes indicator functions to a.s.~indicator functions.

Now suppose $A, B \in \mathcal{F}$ are disjoint.  Then, $\chi_A
\chi_B = 0$.  So, $0 = \int_X \chi_A \chi_B \, d\mu = \langle
\chi_A, \chi_B \rangle = \langle T(\chi_A), T(\chi_B) \rangle =
\int_X T(\chi_A) T(\chi_B) \, d\mu$.  As the integrand is
necessarily nonnegative a.s.[$\mu$], we have $T(\chi_A) T(\chi_B) =
0$ a.s.[$\mu$].

Let $s = \sum_{i=1}^n c_i \chi_{A_i}$ be a simple function, where
$A_j$ are pairwise disjoint.  Then, $s^2 = \sum c_i^2 \chi_{A_i}$.
Now, $T(s) = \sum c_i T(\chi_{A_i})$, which gives $T(s)^2 = \sum_i
\sum_k c_i c_k T(\chi_{A_i}) T(\chi_{A_k}) = \sum_i c_i^2
T(\chi_{A_i})^2 = \sum_i c_i^2 T(\chi_{A_i}) = T(s^2)$ a.s.[$\mu$].

Let $f \in L^2(\mu)$ and $\epsilon > 0$.  Denote by $\|T\|$ the
operator norm of $T$ on $L^1(\mu)$.  Choose a simple function $s$ so
that $|s| \leq |f|$, $\|s - f\|_{L^2(\mu)} <
\epsilon/(4\|f\|_{L^2(\mu)})$ and $\|s^2 - f^2\|_{L^1(\mu)} <
\epsilon/(2 \|T\|)$.  Then,

\begin{equation*}
\begin{split}
\|T(f^2) - &T(f)^2\|_{L^1(\mu)} \\
&\leq \|T(f^2) - T(s^2)\|_{L^1(\mu)} + \|T(s^2) -
T(s)^2\|_{L^1(\mu)}
+ \|T(s)^2 - T(f)^2\|_{L^1(\mu)} \\
&= \|T(f^2 - s^2)\|_{L^1(\mu)} +
\|(T(f) - T(s))(T(f) + T(s)) \|_{L^1(\mu)} \\
&\leq \|T\| \|f^2 - s^2\|_{L^1(\mu)} + \|T(f + s)\|_{L^2(\mu)}
\|T(f - s)\|_{L^2(\mu)} \\
&< \epsilon/2 + 2\|f\|_{L^2(\mu)} \|f - s\|_{L^2(\mu)} \\
&< \epsilon.
\end{split}
\end{equation*}

\noindent As $\epsilon$ is arbitrary, we have $T(f^2) = T(f)^2$
a.s.[$\mu$].
\end{proof}

\vskip0.1in

Our results will focus on a more general sequence. Let $S_n :
\mathcal{M}(X) \rightarrow \mathcal{M}(X)$, $n \in \mathbb{N}$, be a
sequence of linear operators where each $S_n : L^2(\mu) \rightarrow
L^2(\mu)$ is bounded (not necessarily uniformly).  We have the
following Banach principle-type statements.

\begin{thm}\label{thm:bj1} Let $S_n : L^2(\mu) \rightarrow L^2(\mu)$
be bounded.  Suppose that for some $2 \leq p < \infty$, $\sup_n |S_n
f|$ is finite a.s.[$\mu$] for all $f \in L^p(\mu)$. Then, there is a
finite-valued function $\theta(\epsilon)$ so that

\begin{equation*}
\mu \Big\{x \in X : \sup_n |S_n f(x)| > \theta(\epsilon) \Big\} <
\epsilon
\end{equation*}

\noindent for every $\|f\|_{L^p(\mu)} \leq 1$. \end{thm}

\begin{thm}\label{thm:bj2} Let $S_n : L^2(\mu) \rightarrow L^2(\mu)$
be bounded.  Suppose $S_n f$ converges a.s.[$\mu$] for all $f \in
L^\infty(\mu)$.  Then, for each $\epsilon > 0$ and $\eta > 0$, there
exists $\rho(\epsilon, \eta) > 0$ such that

\begin{equation*}
\mu \Big\{ x \in X : \sup_n |S_n f(x)| > \eta \Big\} < \epsilon
\end{equation*}

\noindent for all $\|f\|_{L^\infty(\mu)} \leq 1$ and
$\|f\|_{L^1(\mu)} \leq \rho(\epsilon, \eta)$. \end{thm}

\vskip0.1in

Theorem~\ref{thm:bj1} is a classical result, and
Theorem~\ref{thm:bj2} is due to Bellow and
Jones~\cite{bellowandjones}.  We postpone the proofs of these two
theorems until Section~\ref{sec:banach}. The principal assumption we
make on the sequence $(S_n)$ is that it commutes with a Bourgain
sequence $(T_j)$, that is, $S_n T_j = T_j S_n$ for all $n, j$. More
on this later.

\section{Normal Random Variables}

The proofs of the two entropy results rely heavily on the theory of
Gaussian (or normal) random variables and Gaussian processes.  It is
advantageous at this point to review a few fundamental results.  Let
$(\Omega, \mathcal{B}, P)$ be another probability space.

\begin{defn} We say a random variable $g : \Omega \rightarrow
\mathbb{R}$ is normal (or Gaussian) with mean $m$ and variance
$\sigma^2$ if it has the density function

\begin{equation*}
f(x) = \frac{1}{\sqrt{2\pi \sigma^2}} \exp \left(- \frac{(x -
m)^2}{2\sigma^2} \right),
\end{equation*}

\noindent i.e., $P(g \in A) = \int_A f(x) \, dx$ for all Borel sets
$A$. A normal random variable with mean 0 and variance 1 is called a
standard normal random variable.\end{defn}

Recall, for a random variable $g$ with mean 0, the variance is given
by $\sigma^2 = E(g^2) = \|g\|_{L^2(P)}^2$, where $E(\cdot)$ is
expectation. Also, if $g$ is a normal random variable with mean 0,
then it is centered, that is, $P(g > 0) = P(g < 0) = 1/2$. The
following results are well-known in probability theory, and we
present them without proof.

\begin{lemma}\label{lemma:expect}
If $g$ is a random variable on $\Omega$, then

\begin{equation*}
\int_\Omega |g(\omega)| \, P(d\omega) = \int_0^\infty P(|g| > t) \,
dt \quad \text{and} \quad \int_\Omega |g(\omega)|^2 \, P(d\omega) =
2 \int_0^\infty t P(|g|
> t) \, dt.
\end{equation*}
\end{lemma}

\begin{lemma}\label{lemma:normal1} Let $g$ be a normal random variable
with mean 0.  Then, the moment generating function is given by

\begin{equation*}
\int_\Omega e^{\lambda g(\omega)} \, P(d\omega) = e^{\lambda^2
\sigma^2/2}, \quad \lambda \in \mathbb{R}.
\end{equation*}
\end{lemma}

\begin{lemma}\label{lemma:normal2} For each $1 \leq p < \infty$, there
exists a constant $C_p > 0$ so that $\|g\|_{L^p(P)} \leq C_p
\|g\|_{L^2(P)}$ for all normal random variables $g$ with mean 0.
\end{lemma}

\begin{lemma}\label{lemma:normal3} Let $g_1, g_2, \ldots, g_m$ be
independent standard normal random variables. Then, for any
constants $a_i$, $\sum a_i g_i$ is a normal random variable with
mean 0 and variance $\sum a_i^2$.
\end{lemma}

\vskip0.1in

Lemma~\ref{lemma:expect} is proven by two simple applications of
Fubini's Theorem.  The proof of Lemma~\ref{lemma:normal1} is a
standard result and is found in most probability texts.
Lemma~\ref{lemma:normal3} follows immediately from independence.
Only Lemma~\ref{lemma:normal2} is a somewhat deep result.  In fact,
something stronger is true; the $L^p$ and $L^q$ norms of a normal
random variable are uniformly equivalent for any $1 \leq p, q <
\infty$.  We need only the case $q = 2$.  A proof of the general
result can be found in~\cite{ledoux} (Corollary 3.2).

We now state an important estimate for normal random
variables~\cite{fernique}. The proof is postponed until
Section~\ref{sec:fern}.

\begin{thm}\label{thm:fernique}
Let $G_1, \ldots, G_N$ be normal random variables each with mean 0.
If for some constant $s$ we have $P\{ \omega \in \Omega : \sup_n
|G_n(\omega)| \leq s \} \geq 1/2$, then $\|\sup_n |G_n|\|_{L^1(P)}
\leq 6s$.
\end{thm}

\vskip0.1in

We now turn our attention to Gaussian processes.

\begin{defn} Let $T$ be a countable indexing set.  We say a
collection of random variables $(G_t : t \in T)$ is a Gaussian
process if each $G_t$ has mean 0 and all finite linear combinations
$\sum_t a_t G_t$ are normal random variables.
\end{defn}

Note that this definition is not entirely standard, in particular
the requirement that each $G_t$ have mean 0. It is added here
because throughout all Gaussian processes we deal with have this
property, and it makes life simpler later on.

If each $G_t$ is itself a finite linear combination of mean 0 normal
random variables, then $(G_t : t \in T)$ is trivially a Gaussian
process. Define a pseudo-metric on $T$ by $d_G(s,t) = \|G_s -
G_t\|_{L^2(P)}$. Denote the entropy number of $T$ by
$N(T,d_G,\delta)$.  The following fundamental result is Sudakov's
inequality.

\begin{thm}
There exists a universal constant $R$ such that if $(G_t : t \in T)$
is a Gaussian process, then

\begin{equation*}
\sup_{\delta > 0} \, \delta \sqrt{ \log N(T,d_G,\delta)} \leq R
\left\| \sup_{t \in T} |G_t| \right\|_{L^1(P)}.
\end{equation*}
\end{thm}

\vskip0.1in

For the remainder of the paper, the notation $R$ is fixed on this
constant.  A proof of Sudakov's inequality can be found
in~\cite{ledoux} (Theorem 3.18).

\section{The First Entropy Result}

Recall, we consider a sequence $S_n : L^2(\mu) \rightarrow L^2(\mu)$
of bounded operators.  For $f \in L^2(\mu)$, define a pseudo-metric
on $\mathbb{N}$ by $d_f(n,n') = \|S_n f - S_{n'} f\|_{L^2(\mu)}$.
For $\delta > 0$, let $N_f(\delta) := N(\mathbb{N}, d_f, \delta)$,
that is, the $\delta$-entropy number of the set $\{S_n f : n \in
\mathbb{N}\}$ in $L^2(\mu)$. We now state and prove the first of
Bourgain's entropy results.

\begin{prop} Let $S_n : L^2(\mu) \rightarrow L^2(\mu)$ be bounded
(not necessarily uniformly), and assume $(S_n)$ commutes with a
Bourgain sequence $(T_j)$. Suppose that for some $1 \leq p <
\infty$, $\sup_n |S_n f| < \infty$ a.s.[$\mu$] for all $f \in
L^p(\mu)$. Then, there exists a constant $C > 0$ such that $\delta
(\log N_f(\delta))^{1/2} \leq C \|f\|_{L^2(\mu)}$ for all $\delta >
0$ and $f \in L^2(\mu)$.  \end{prop}

\begin{proof} As $(X, \mu)$ is a probability space, $L^p(\mu)
\supset L^q(\mu)$ when $p < q$.  So, if $p < 2$, then $\sup_n |S_n
f| < \infty$ a.s.[$\mu$] for all $f \in L^p(\mu) \supset L^q(\mu)$
for any $q \geq 2$.  Therefore, assume without loss of generality
that $p \geq 2$.

For $M \in \mathbb{N}$, let $\overline{M} = \{1, 2, \ldots, M\}$.
Note, $\sup_{M} N(\overline{M}, d_f, \delta) = N_f(\delta)$.
Therefore, it suffices to find $C$, independent of $M$, such that
$(\log N(\overline{M},d_f,\delta))^{1/2} \leq C\|f\|_{L^2(\mu)}$ for
all $f, \delta$.  Fix $M \in \mathbb{N}$.

Suppose we could show $\delta \log(N(\overline{M}, d_f,
\delta))^{1/2} \leq C\|f\|_{L^2(\mu)}$ for all $f \in L^\infty(\mu)$
and all $\delta > 0$.  Let $f \in L^2(\mu)$ and $\delta > 0$.  Let
$D = \max (\|S_1\|, \ldots, \|S_M\|)$ be the maximum of the
$L^2(\mu)$ operator norms.  Choose $f_1 \in L^\infty(\mu)$ with
$|f_1| \leq |f|$ and $\|f - f_1\|_{L^2(\mu)} < \delta/(2D)$. Then,
$\|S_n f - S_n f_1\|_{L^2(\mu)} < \delta/2$ for all $n \in
\overline{M}$ and $N(\overline{M},d_f,\delta) \leq
N(\overline{M},d_{f_1},\delta/2)$. Hence, $\delta(\log
N(\overline{M},d_f,\delta))^{1/2} \leq 2C\|f_1\|_{L^2(\mu)} \leq
2C\|f\|_{L^2(\mu)}$, and we have the desired estimate with $2C$.
Therefore, it suffices to prove the result for all $L^\infty(\mu)$
functions.  Fix $f \in L^\infty(\mu)$.

We will fix $J$ at some large integer.  By the mean ergodic
condition on $T_j$, we have

\begin{equation*}
J^{-1} \sum_{j=1}^J T_j(f^2) \rightarrow \|f^2\|_{L^1(\mu)} =
\|f\|_{L^2(\mu)}^2
\end{equation*}

\noindent in $L^1(\mu)$-norm. But, $|J^{-1} \sum T_j (f^2)| \leq
\|f\|_{L^\infty(\mu)}^2$ a.s.[$\mu$] by Lemma~\ref{lemma:T}.  Of
course, if a sequence of functions $h_n$ converges to a constant $c$
in $L^1(\mu)$-norm and $|h_n| \leq B$ a.s.[$\mu$] for all $n$, it
follows $h_n \rightarrow c$ in $L^q(\mu)$-norm for all $1 \leq q <
\infty$. Hence, as $p \geq 2$, choose $J$ so large that

\begin{equation}\label{eq:bigJ1}
\left\| \frac{1}{J} \sum_{j=1}^J T_j(f^2) \right\|_{L^{p/2}(\mu)}
\leq 2\|f\|_{L^2(\mu)}^2.
\end{equation}

Similarly, for each pair $n, n' \in \overline{M}$, $J^{-1} \sum
T_j(S_n f - S_{n'} f)^2 \rightarrow \|S_n f - S_{n'}
f\|_{L^2(\mu)}^2 = d_f(n,n')^2$ in $L^1(\mu)$-norm, and thus in
probability. So, for $J$ large enough,

\begin{equation*}
\mu(\mathcal{C}_{n,n'}) := \mu \left\{x \in X : J^{-1} \sum_{j=1}^J
T_j (S_n f - S_{n'} f)^2(x) \geq \frac{1}{4} d_f(n,n')^2 \right\} >
1 - \frac{1}{16M^2}.
\end{equation*}

\noindent Choose $J$ big enough so that this holds for each pair $n,
n'$.

Let $g_1, \ldots, g_J$ be a sequence of independent standard normal
random variables on a probability space $(\Omega, \mathcal{B}, P)$.
Define the functions $F, F^*$ on the product space $X \times \Omega$
by

\begin{equation*}
F(x,\omega) = J^{-1/2} \sum_{j=1}^J g_j(\omega) T_j f(x) \qquad
\text{and} \qquad F^*(x,\omega) = \sup_{n \in \overline{M}} |S_n
F(x,\omega)|.
\end{equation*}

\noindent By the commutativity assumption, $S_n F(x,\omega) =
J^{-1/2} \sum g_j(\omega) T_j S_n f(x)$.  Note, for each fixed $x$
such that $T_j S_n f(x)$ is finite for all $j, n$, $(S_n F(x,\cdot)
: n \in \overline{M})$ is a Gaussian process. The focus of the proof
will be finding the ``correct" $x$ to fix.

We define four sets $\mathcal{A}, \mathcal{B}, \mathcal{C},
\mathcal{D} \subset X$. First, let

\begin{equation*}
\mathcal{A} = \Big\{x \in X : |T_j S_n f(x)| < \infty \, \text{for
all } 1 \leq j \leq J, \, n \in \overline{M} \Big\}.
\end{equation*}

\noindent As each $T_j S_n f \in L^2(\mu)$, it is clear that
$\mu(\mathcal{A}) = 1$.  Set

\begin{equation*}
\mathcal{B} = \Big\{ x \in X : T_j(S_n f - S_{n'}f)^2(x) = T_j(S_n f
- S_{n'}f)(x)^2 \,\, \text{for all } 1 \leq j \leq J, \, n, n' \in
\overline{M} \Big\}
\end{equation*}

\noindent By Lemma~\ref{lemma:T}, $\mu(\mathcal{B}) = 1$. Let
$\mathcal{C} = \bigcap_{n, n'} \mathcal{C}_{n,n'}$, i.e.,

\begin{equation*}
\mathcal{C} = \left\{ x \in X : \left( J^{-1} \sum_{j=1}^J T_j(S_n f
- S_{n'} f)^2(x) \right)^{1/2} \geq \frac{1}{2} d_f(n,n') \text{ for
all } n, n' \in \overline{M} \right\}.
\end{equation*}

\noindent Now,

\begin{equation*}
\begin{split}
\mu \bigg( \bigcup_{n,n' \in \overline{M}} \mathcal{C}^c_{n,n'}
\bigg) &\leq \sum_{n,n'} \mu(\mathcal{C}^c_{n,n'}) < \sum_{n,n'}
\frac{1}{16M^2} = \frac{1}{16}.
\end{split}
\end{equation*}

\noindent Equivalently, $\mu(\mathcal{C}) > 15/16$.

Finally, let $R' = 5\sqrt{2} C_p \, \theta(1/4)$, where where
$\theta$ is the finite-valued function from Theorem~\ref{thm:bj1}
and $C_p$ is the constant from Lemma~\ref{lemma:normal2}.  Note,
$R'$ does not depend on $f$ or $M$. Set

\begin{equation*}
\mathcal{D} = \Big\{ x \in X : P\big\{ \omega \in \Omega :
F^*(x,\omega) \leq R'\|f\|_{L^2(\mu)} \big\} \geq 1/2 \Big\}.
\end{equation*}

Now, by Lemmas~\ref{lemma:T}, \ref{lemma:normal2}, and
\ref{lemma:normal3}, and (\ref{eq:bigJ1}), we have

\begin{equation*}
\begin{split}
\int_\Omega \|F(\cdot,\omega)\|_{L^p(\mu)} \, P(d\omega) &=
\int_\Omega \left( \int_X |F(x,\omega)|^p \, \mu(dx) \right)^{1/p}
\, P(d\omega) \\
&\leq \left( \int_\Omega \int_X |F(x,\omega)|^p \, \mu(dx) \,
P(d\omega) \right)^{1/p} \\
&= \left( \int_X \int_\Omega \left| \sum_{j=1}^J J^{-1/2}
g_j(\omega) T_j f(x) \right|^p \, P(d\omega) \, \mu(dx) \right)^{1/p} \\
&= \left( \int_X \left\| \sum_{j=1}^J J^{-1/2} g_j(\cdot) T_j f(x)
\right\|_{L^p(P)}^p \, \mu(dx) \right)^{1/p} \\
&\leq \left( \int_X C_p^p \left\| \sum_{j=1}^J J^{-1/2} g_j(\cdot)
T_j f(x) \right\|_{L^2(P)}^p \, \mu(dx) \right)^{1/p} \\
&= C_p \left( \int_X \left| \frac{1}{J} \sum_{j=1}^J (T_j f)(x)^2
\right|^{p/2} \, \mu(dx) \right)^{1/p} \\
&= C_p \left\| \frac{1}{J} \sum_{j=1}^J T_j(f^2)
\right\|_{L^{p/2}(\mu)}^{1/2} \\
&\leq C_p \sqrt{2} \|f\|_{L^2(\mu)}.
\end{split}
\end{equation*}

\noindent By Chebyshev's inequality, $P \{\omega :
\|F(\cdot,\omega)\|_{L^p(\mu)} > 5\sqrt{2}C_p\|f\|_{L^2(\mu)} \}
\leq 1/5 < 1/4$, or equivalently,

\begin{equation*}
P \Big\{\omega \in \Omega : \|F(\cdot,\omega)\|_{L^p(\mu)} \leq
5\sqrt{2}C_p \|f\|_{L^2(\mu)}\Big\} > 3/4.
\end{equation*}

\noindent Fix an $\omega$ in the above set.  Then, $F(\cdot,\omega)
\in L^p(\mu)$.  By Theorem~\ref{thm:bj1},

\begin{equation*}
\mu \Big\{x \in X : \sup_{n \in \mathbb{N}} |S_n F(x,\omega)| \leq
\|F(\cdot,\omega)\|_{L^p(\mu)} \theta(1/4) \Big\} > 3/4.
\end{equation*}

\noindent Of course, $F^*(x,\omega) \leq \sup_n |S_n F(x,\omega)|$.
Further, we have a bound on $\|F(\cdot,\omega)\|_{L^p(\mu)}$ by the
choice of $\omega$. Thus,

\begin{equation*}
\mu \Big\{x \in X : F^*(x,\omega) \leq 5\sqrt{2}C_p \,
\theta(1/4)\|f\|_{L^2(\mu)} \Big\} > 3/4.
\end{equation*}

\noindent As this holds for all such $\omega$, we have

\begin{equation*}
P \Big\{ \omega \in \Omega : \mu \big\{x \in X : F^*(x,\omega) \leq
R'\|f\|_{L^2(\mu)} \big\} > 3/4 \Big\} > 3/4.
\end{equation*}

We now apply Fubini's theorem.  In particular,

\begin{equation*}
\begin{split}
3/4 &< \int_\Omega \chi_{\{\mu(F^*(x,\omega) \leq
R'\|f\|_{L^2(\mu)})
>3/4\}}(\omega) \, P(d\omega) \leq \int_\Omega \frac{4}{3} \mu\{x :
F^*(x,\omega) \leq R'\|f\|_{L^2(\mu)}\} \, P(d\omega) \\
&= \frac{4}{3} \int_\Omega \int_X \chi_{\{F^*(x,\omega) \leq
R'\|f\|_{L^2(\mu)}\}}(x,\omega) \, \mu(dx) \, P(d\omega) \\
&= \frac{4}{3} \int_X \int_\Omega \chi_{\{F^*(x,\omega) \leq
R'\|f\|_{L^2(\mu)}\}}(x,\omega) \, P(d\omega)
\, \mu(dx) \\
&= \frac{4}{3} \int_X P\{\omega : F^*(x,\omega) \leq
R'\|f\|_{L^2(\mu)}\} \, \mu(dx) \\
&= \frac{4}{3} \bigg[ \int_\mathcal{D} P\{F^*(x,\omega) \leq
R'\|f\|_{L^2(\mu)}\} \, \mu(dx) + \int_{\mathcal{D}^c}
P\{F^*(x,\omega)
\leq R'\|f\|_{L^2(\mu)}\} \, \mu(dx) \bigg] \\
&\leq \frac{4}{3}(\mu(\mathcal{D}) + 1/2).
\end{split}
\end{equation*}

\noindent The last line follows from the definition of
$\mathcal{D}$. This gives $\mu(\mathcal{D}) > 1/16$.

The estimates $\mu(\mathcal{A}) = 1$, $\mu(\mathcal{B}) = 1$,
$\mu(\mathcal{C}) > 15/16$, and $\mu(\mathcal{D}) > 1/16$ together
imply that $\mu(\mathcal{A} \cap \mathcal{B} \cap \mathcal{C} \cap
\mathcal{D}) > 0$. Fix $\overline{x} \in \mathcal{A} \cap
\mathcal{B} \cap \mathcal{C} \cap \mathcal{D}$. Define $G_n(\omega)
= S_n F(\overline{x},\omega) = J^{-1/2} \sum g_j(\omega) T_j S_n
f(\overline{x})$.  As $\overline{x} \in \mathcal{A}$, $(G_n : n \in
\overline{M})$ is a Gaussian process.  By Sudakov's inequality and
Theorem~\ref{thm:fernique}, and because $\overline{x} \in
\mathcal{D}$, we see that

\begin{equation*}
\begin{split}
\sup_{\delta > 0} \delta (\log N(\overline{M},d_G,\delta))^{1/2}
&\leq R \int_\Omega \sup_{n \in \overline{M}} |G_n(\omega)| \,
P(d\omega) \\
&= R \int_\Omega F^*(\overline{x},\omega) \, P(d\omega) \leq
6RR'\|f\|_{L^2(\mu)}.
\end{split}
\end{equation*}

\noindent On the other hand, each $G_n - G_{n'}$ is a linear
combination of independent standard normal random variables. It
follows from Lemma~\ref{lemma:normal3} again and because
$\overline{x} \in \mathcal{B} \cap \mathcal{C}$ that

\begin{equation*}
\begin{split}
d_G(n,n') &= \|G_n - G_{n'}\|_{L^2(P)} = \bigg( \frac{1}{J}
\sum_{j=1}^J \Big(T_j S_n f(\overline{x}) - T_j S_{n'}
f(\overline{x})
\Big)^2 \bigg)^{1/2} \\
&= \left( \frac{1}{J} \sum_{j=1}^J T_j (S_n f - S_{n'}
f)^2(\overline{x}) \right)^{1/2} \\
&\geq \frac{1}{2} d_f(n,n').
\end{split}
\end{equation*}

\noindent This implies $N(\overline{M},d_f,\delta) \leq
N(\overline{M},d_G,\delta/2)$ for all $\delta > 0$.  Hence, $\delta
(\log N(\overline{M},d_f,\delta))^{1/2} \leq 12RR'\|f\|_{L^2(\mu)}$.
We note that $12RR'$ is universal, and does not depend on $M$ or
$f$. As $f \in L^\infty(\mu)$ was arbitrary, this holds all
a.s.~bounded functions.  By our earlier note, $\delta (\log
N(\overline{M},d_f,\delta))^{1/2} \leq 24RR'\|f\|_{L^2(\mu)} =:
C\|f\|_{L^2(\mu)}$ for all $f \in L^2(\mu)$ and all $\delta > 0$.
Taking the supremum over $M$, $\delta (\log N_f(\delta))^{1/2} \leq
C\|f\|_{L^2(\mu)}$.
\end{proof}

\section{The Second Entropy Result}

It will now be necessary to assume the $(S_n)$ are uniformly
bounded.  Of course, by dividing out a constant, we may assume each
$S_n$ is an $L^2(\mu)$-contraction, i.e., $\|S_n f\|_{L^2(\mu)} \leq
\|f\|_{L^2(\mu)}$ for all $n$.

\begin{prop} Let $S_n$ be a sequence of $L^2(\mu)$ contractions that
commute with a Bourgain sequence $(T_j)$.  Suppose $S_n f$ converges
a.s.[$\mu$] for all $f \in L^\infty(\mu)$.  Then, there exists a
finite-valued function $C(\delta)$ such that $N_f(\delta) \leq
C(\delta)$ for all $\delta > 0$ and $\|f\|_{L^2(\mu)} \leq 1$.
\end{prop}

\begin{proof}
Suppose we can show the uniform entropy estimate for all $f \in
L^\infty(\mu)$, $\|f\|_{L^2(\mu)} \leq 1$.  Let $\|f\|_{L^2(\mu)}
\leq 1$ and $\delta > 0$. Choose $f_1 \in L^\infty(\mu)$, $|f_1|
\leq |f|$ such that $\|f - f_1\|_{L^2(\mu)} < \delta/2$. Then,
$\|S_n f - S_n f_1\|_{L^2(\mu)} < \delta/2$ for all $n$ and
$N_f(\delta) \leq N_{f_1}(\delta/2) \leq C(\delta/2)$. As $\delta$
and $f$ are arbitrary, we have the uniform estimate with the
function $C_0(\delta) = C(\delta/2)$. It therefore suffices to prove
the result for a.s.~bounded functions.

We proceed by contradiction.  Suppose not, i.e., suppose there is
some $\delta > 0$ such that $N_f(\delta)$ is unbounded over all such
$f$. Define the constant $R' = \frac{\delta}{25R}$. As per
Theorem~\ref{thm:bj2}, pick the constant $\rho(1/10, R'/10)$. Choose
$K \in \mathbb{N}$, $K > 1$ big enough so that $\frac{1}{6}(R' -
1000(\log K)^{-1/2}) > R'/10$ and $2(\log K)^{-1/2} <
\rho(1/10,R'/10)$.

Now, by our assumption, there is some $f \in L^\infty$,
$\|f\|_{L^2(\mu)} \leq 1$ such that $N_f(\delta) > K$. In
particular, there is a subset $I \subset \mathbb{N}$ with $|I| = K$
(cardinality) and $\|S_n f - S_{n'} f\|_{L^2(\mu)} > \delta$ for all
$n \not= n' \in I$.

As before, we will need to choose an appropriately large $J$. First,
denote $B = \|f\|_{L^\infty(\mu)}$. Fix a number $T > 0$ such that
$T > 3 \sqrt{ \log K}$ and $\exp(\frac{-T^2}{2B^2}) \leq
\frac{7}{2KB^2}$. Note, the quantity $e^{\lambda^2 (2-B^2)} -
e^{\lambda^2(1-B^2)}$ is strictly positive for all $\lambda \in
[\sqrt{\log K}, T/3]$. Let $\gamma > 0$ be the minimum value of this
quantity for $\lambda$ in this interval. As $J^{-1} \sum T_j (f^2)
\rightarrow \|f^2\|_{L^1(\mu)} = \|f\|_{L^2(\mu)}^2$ in
$L^{1}(\mu)$-norm, it converges in probability.  So, pick $J$ large
enough so that

\begin{equation}\label{eq:bigJ3}
\mu(Y) := \mu \left\{x \in X : \left| \frac{1}{J} \sum_{j=1}^J
T_j(f^2)(x) - \|f\|_{L^2(\mu)}^2 \right| > 1 \right\} < \gamma.
\end{equation}

Recall from the proof of Proposition 1, $J^{-1} \sum T_j(S_nf -
S_{n'}f)^2 \rightarrow \|S_n f - S_{n'}f\|_{L^2(\mu)}^2$ in
probability for each pair $n, n' \in I$. So, just as we did before,
take $J$ big enough so that if

\begin{gather}\label{eq:bigJ4}
Z_1 = \left\{ x \in X : \left( J^{-1} \sum_{j=1}^J T_j(S_n f -
S_{n'} f)^2(x) \right)^{1/2} \geq \frac{1}{2} \|S_n f - S_{n'}
f\|_{L^2(\mu)} \text{ for all } n, n' \in I \right\}\nonumber \\
\text{then } \,\, \mu(Z_1) > 4/5.
\end{gather}

\vskip0.1in

Again, define $F(x,\omega) = J^{-1/2} \sum_{j=1}^J g_j(\omega) T_j
f(x)$. Write $F(x,\omega) = \varphi(x,\omega) + H(x,\omega)$ where
\begin{equation*}
\varphi(x,\omega) = F(x,\omega) \chi_{\{|F(x,\omega)| \leq 6
\sqrt{\log K}\}} \quad \text{and} \quad H(x,\omega) = F(x,\omega)
\chi_{\{|F(x,\omega)| > 6 \sqrt{\log K}\}}.
\end{equation*}

\noindent Define three subsets of $\Omega$ by

\begin{gather*}
\mathcal{A} = \bigg\{\omega \in \Omega : \int_X \sup_{n \in I} |S_n
H(x,\omega)| \,
\mu(dx) \leq 90 \bigg\}, \\
\mathcal{B} = \bigg\{\omega \in \Omega :  \mu \big\{x \in X :
\sup_{n \in I} |S_n F(x,\omega)| > R'(\log K)^{1/2} \big\} > 1/5 \bigg\}, \\
\mathcal{C} = \bigg\{\omega \in \Omega :  \int_X |\varphi(x,\omega)|
\, \mu(dx) \leq 12 \bigg\}.
\end{gather*}

Suppose for the moment that we could choose $\overline{\omega} \in
\mathcal{A} \, \cap \, \mathcal{B} \, \cap \, \mathcal{C}$. Define
$\psi(x) = \frac{1}{6} (\log K)^{-1/2}
\varphi(x,\overline{\omega})$. Simply from the definition of
$\varphi$, it follows

\begin{equation*}
|\psi| \leq 1.
\end{equation*}

\noindent As $\overline{\omega} \in \mathcal{C}$,

\begin{equation*}
\int_X |\psi(x)| \, \mu(dx) \leq 2(\log K)^{-1/2} < \rho(1/10,
R'/10).
\end{equation*}

\noindent By Theorem~\ref{thm:bj2}, we have that

\begin{equation}\label{eq:psi1}
\mu \left\{ x \in X : \sup_n |S_n \psi(x)| > R'/10 \right\} < 1/10.
\end{equation}

On the other hand, as $\overline{\omega} \in \mathcal{A}$, we have
by Chebyshev that

\begin{equation*}
\mu \left\{x \in X : \sup_{n \in I} |S_n H(x,\overline{\omega})| >
1000 \right\} \leq 90/1000 < 1/10,
\end{equation*}

\noindent or equivalently

\begin{equation}\label{eq:cgame1}
\mu \left\{x \in X : \sup_{n \in I} |S_n H(x,\overline{\omega})|
\leq 1000 \right\} > 9/10.
\end{equation}

\noindent As $\overline{\omega} \in \mathcal{B}$,

\begin{equation}\label{eq:cgame2}
\mu \left\{ x \in X : \sup_{n \in I} |S_n F(x,\overline{\omega})|
> R'(\log K)^{1/2} \right\} > 1/5.
\end{equation}

\noindent Now, $\sup_n |S_n \varphi| \geq \sup_I |S_n \varphi| \geq
\sup_I |S_n F| - \sup_I |S_n H|$.  So, taking the intersection of
the sets in (\ref{eq:cgame1}) and (\ref{eq:cgame2}), we have

\begin{equation*}
\mu \left\{x \in X : \sup_n |S_n \varphi(x,\overline{\omega})| \geq
R'(\log K)^{1/2} - 1000 \right\}
> 1/10.
\end{equation*}

\noindent Applying the definition of $\psi$ and the choice of $K$,

\begin{equation*}
\mu \left\{ x \in X : \sup_n |S_n \psi(x)| > R'/10 \right\} > 1/10.
\end{equation*}

\noindent This clearly contradicts (\ref{eq:psi1}).  Therefore, it
suffices to find such an $\overline{\omega}$.

\vskip0.1in

\noindent {\bf Estimate of $\mathcal{A}$}.

\vskip0.1in

Fix $\lambda \in [\log K, T/3]$.  By considering $F(x,\omega)$ as a
normal random variable (in $\omega$) with variance $J^{-1} \sum
T_jf(x)^2$, it follows from Lemma~\ref{lemma:normal1} that

\begin{equation*}
\begin{split}
\int_X \int_\Omega \exp(\lambda F(x,\omega)) \, P(d\omega) \,
\mu(dx) &= \int_X \exp \left( \frac{\lambda^2}{2 J} \sum_{j=1}^J
T_jf(x)^2 \right) \, \mu(dx) \\
&= \int_X \exp \left( \frac{\lambda^2}{2 J} \sum_{j=1}^J T_j(f^2)(x)
\right) \, \mu(dx).
\end{split}
\end{equation*}

\noindent Now, by (\ref{eq:bigJ3}), we see $J^{-1} \sum T_j(f^2)
\leq \|f\|_{L^2(\mu)}^2 + 1 \leq 2$ on $Y^c$. On the other hand, we
have $J^{-1} \sum T_j(f^2) \leq B^2$ a.s.[$\mu$].  So,

\begin{equation*}
\begin{split}
\int_X \int_\Omega \exp(\lambda F(x,\omega)) \, P(d\omega) \,
\mu(dx) &= \int_X \exp \left( \frac{\lambda^2}{2 J} \sum_{j=1}^J
T_j(f^2)(x)
\right) \, \mu(dx) \\
&\leq \int_{Y^c} e^{\lambda^2} \, \mu(dx) + \int_Y e^{\lambda^2
B^2} \, \mu(dx) \\
&\leq e^{\lambda^2} + \gamma e^{\lambda^2 B^2} \\
&\leq e^{\lambda^2} + (e^{\lambda^2 (2-B^2)} - e^{\lambda^2(1-B^2)})
e^{\lambda^2 B^2} = e^{2\lambda^2}.
\end{split}
\end{equation*}

\noindent Define $\mu_t(\omega) = \mu \{x \in X : |F(x,\omega)| > t
\}$.  Then, for all $t > 0$, we have from above that

\begin{equation*}
\begin{split}
e^{\lambda t} \int_\Omega \mu_t(\omega) \, P(d\omega) & =
e^{\lambda t} \int_X P\{ |F(x,\omega)| > t\} \, \mu(dx) \\
&= 2 e^{\lambda t} \int_X P\{F(x,\omega) > t\} \, \mu(dx) \\
&= 2 \int_X \int_{\{F(x,\omega) > t\}} e^{\lambda t} \,
P(d\omega) \, \mu(dx) \\
&\leq 2 \int_X \int_\Omega e^{\lambda F(x,\omega)} \, P(d\omega)
\, \mu(dx) \\
&\leq 2 e^{2\lambda^2}.
\end{split}
\end{equation*}

\noindent Set $t = 3\lambda$.  Then, $\int_\Omega \mu_t(\omega) \,
P(d\omega) \leq 2e^{-t^2/9}$, and this holds for all $t \in [3
\sqrt{\log K}, T]$.

On the other hand, $J^{-1} \sum T_j(f^2) \leq B^2$ a.s.[$\mu$]. So,
for all $t, \lambda > 0$,

\begin{equation*}
\begin{split}
e^{\lambda t} \int_\Omega \mu_t(\omega) \, P(d\omega) &\leq 2 \int_X
\int_\Omega e^{\lambda F(x,\omega)} \, P(d\omega) \mu(dt)
\\
&= 2 \int_X \exp \left( \frac{\lambda^2}{2J} \sum_{j=1}^J
T_j(f^2)(x) \right) \,
\mu(dx) \\
&\leq 2 \int_X \exp \left(\frac{\lambda^2 B^2}{2} \right) \, \mu(dx) \\
&= 2 e^{\lambda^2 B^2 / 2}.
\end{split}
\end{equation*}

\noindent Set $t = \lambda B^2$ to see $\int_X \mu_t(\omega) \,
P(d\omega) \leq 2e^{-\lambda^2 B^2/2} = 2e^{-t^2/(2B^2)}$ for all $t
> 0$.

From Lemma~\ref{lemma:expect},

\begin{equation}\label{eq:mess1}
\begin{split}
\int_\Omega \int_X |H(x,\omega)|^2 \, \mu(dx) \, P(d\omega) &= 2
\int_\Omega \int_0^\infty t \mu\{|H(x,\omega)| > t\} \, dt \,
P(d\omega) \\
&= 2 \int_\Omega \int_0^{3 \sqrt{\log K}} t \mu\{|H(x,\omega)|
> t\} \, dt \, P(d\omega) \\
&\quad + 2 \int_\Omega \int_{3 \sqrt{\log K}}^\infty
t \mu\{|H(x,\omega)| > t\} \, dt \, P(d\omega) \\
\end{split}
\end{equation}

\noindent By definition of $H$, and an application of Chebyshev,

\begin{equation}\label{eq:mess2}
\begin{split}
2 \int_\Omega \int_0^{3 \sqrt{\log K}} &t \mu\{|H(x,\omega)|
> t\} \, dt \, P(d\omega) \\
&= 2 \int_\Omega \int_0^{3 \sqrt{\log
K}} t \mu \big\{|H(x,\omega)| > 6 \sqrt{\log K} \big\} \, dt P(d\omega) \\
&\leq 2\int_\Omega \int_0^{3 \sqrt{\log K}} t \left( \frac{1}{6
\sqrt{\log K}} \right)^2 \left( \int_X |H(x,\omega)|^2 \, \mu(dx)
\right) dt \, P(d\omega) \\
&= \frac{1}{4} \int_\Omega \int_X |H(x,\omega)|^2 \, \mu(dx) \,
P(d\omega) \\
&\leq \frac{1}{2} \int_\Omega \int_X |H(x,\omega)|^2 \, \mu(dx) \,
P(d\omega).
\end{split}
\end{equation}

\noindent As $|H| \leq |F|$ everywhere, $\mu\{|H(x,\omega)| > t\}
\leq \mu_t(\omega)$ for all $t$ and $\omega$.  So,

\begin{equation}\label{eq:mess3}
\begin{split}
2 \int_\Omega \int_{3 \sqrt{\log K}}^\infty &t \mu\{|H(x,\omega)|
> t\} \, dt \, P(d\omega) \\
&\leq 2 \int_\Omega \int_{3 \sqrt{\log K}}^\infty t \mu_t(\omega)
\, dt \, P(d\omega) \\
&= 2 \int_\Omega \int_{3 \sqrt{\log K}}^T t \mu_t(\omega) \, dt \,
P(d\omega) + 2 \int_\Omega \int_T^\infty t \mu_t(\omega) \, dt \,
P(d\omega) \\
&\leq 4 \int_{3 \sqrt{\log K}}^T t e^{-t^2/9} \, dt + 4
\int_T^\infty t e^{-t^2/(2B^2)} \, dt \\
&= 18K^{-1} - \frac{9}{2} e^{-T^2/9} + 4B^2e^{-T^2/(2B^2)} \\
&\leq 32K^{-1}
\end{split}
\end{equation}

\noindent by the choice of $T$.  Combining (\ref{eq:mess1}),
(\ref{eq:mess2}), and (\ref{eq:mess3}),

\begin{equation*}
\frac{1}{2} \int_\Omega \int_X |H(x,\omega)|^2 \, \mu(dx) \,
P(d\omega) \leq 32K^{-1}.
\end{equation*}

\noindent Stated another way, $\int_\Omega
\|H(\cdot,\omega)\|_{L^2(\mu)}^2 \, P(d\omega) \leq 64/K$.  This
gives

\begin{equation*}
\begin{split}
\int_\Omega \int_X \sup_{n \in I} |S_n H(x,\omega)| \, \mu(dx) \,
P(d\omega) &= \int_\Omega \Big\| \sup_{n \in I} |S_n
H(\cdot,\omega)| \Big\|_{L^1(\mu)} \, P(d\omega) \\
&\leq \int_\Omega \Big\| \sup_{n \in I} |S_n
H(\cdot,\omega)| \Big\|_{L^2(\mu)} \, P(d\omega)\\
&= \int_\Omega \Big( \int_X (\sup_{n \in I} |S_n
H(x,\omega)|)^2 \, \mu(dx) \Big)^{1/2} \, P(d\omega)\\
&\leq \int_\Omega \Big( \int_X \sum_{n \in I} |S_n
H(x,\omega)|^2 \, \mu(dx) \Big)^{1/2} \, P(d\omega)\\
&= \int_\Omega \Big( \sum_{n \in I} \|S_n
H(\cdot,\omega)\|_{L^2(\mu)}^2 \Big)^{1/2} \, P(d\omega) \\
&\leq \int_\Omega \Big( \sum_{n \in I} \| H(\cdot,\omega)
\|_{L^2(\mu)}^2 \Big)^{1/2} \, P(d\omega) \\
&= K^{1/2} \int_\Omega \|H(\cdot,\omega)\|_{L^2(\mu)} \,
P(d\omega) \\
&\leq K^{1/2} \left( \int_\Omega \|H(\cdot,\omega)\|_{L^2(\mu)}^2
\, P(d\omega) \right)^{1/2} \\
&\leq K^{1/2} (64/K)^{1/2} = 8.
\end{split}
\end{equation*}

\noindent It follows by Chebyshev that $P \{\int_X \sup_{n \in I}
|S_n H(x,\omega)| \, \mu(dx) > 90\} \leq 8/90 < 1/10$, or
equivalently, $P(\mathcal{A}) > 9/10$.

\vskip0.1in

\noindent {\bf Estimate of $\mathcal{B}$}.

\vskip0.1in

Denote $F^*(x,\omega) = \sup_{n \in I} |S_n F(x,\omega)|$.  Define
the sets $Z_2, Z_3 \subset X$ by

\begin{gather*}
Z_2 = \Big\{ x \in X : |T_j S_n f(x)| < \infty \,\, \text{for all }
1 \leq j \leq J, \,\, n \in I \Big\}, \\
Z_3 = \Big\{ x \in X : T_j(S_n f - S_{n'}f)^2(x) = T_j(S_n f -
S_{n'}f)(x)^2 \,\, \text{for all } 1 \leq j \leq J, \, n, n' \in I
\Big\}.
\end{gather*}

\noindent As in the proof of Proposition 1, $\mu(Z_2) = \mu(Z_3) =
1$. Let $Z = Z_1 \cap Z_2 \cap Z_3$, so that $\mu(Z) > 4/5$ by
(\ref{eq:bigJ4}).

Fix $x \in Z$.  Define a Gaussian process by $G_n(\omega) = S_n
F(x,\omega)$ and let $d_G(n,n') = \|G_n - G_{n'}\|_{L^2(P)}$ as
before. By the definition of $Z$ and the original hypothesis,

\begin{equation*}
\begin{split}
d_G(n,n') &= \left( J^{-1} \sum_{j=1}^J T_j(S_n f - S_{n'} f)(x)^2
\right)^{1/2} \\
&= \left( J^{-1} \sum_{j=1}^J T_j(S_n f - S_{n'} f)^2(x)
\right)^{1/2} \geq \frac{1}{2} \|S_n f - S_{n'} f\|_{L^2(\mu)} >
\delta/2
\end{split}
\end{equation*}

\noindent for all $n \not= n' \in I$. So, $N(I,d_G,\delta/4) = |I| =
K$. By Sudakov's inequality,

\begin{equation*}
\begin{split}
\int_\Omega F^*(x,\omega) \, P(d\omega) &\geq \frac{1}{R}
\frac{\delta}{4} \big[ \log N(I,d_G,\delta/4) \big]^{1/2} \\
&= 6 \left( \frac{\delta}{24 R} \right) (\log K)^{1/2} > 6R' (\log
K)^{1/2}.
\end{split}
\end{equation*}

\noindent It follows from Theorem~\ref{thm:fernique} that $P\{\omega
: F^*(x,\omega) \leq R' (\log K)^{1/2}\} < 1/2$, or equivalently,
$P\{\omega : F^*(x,\omega) > R' (\log K)^{1/2}\} > 1/2$.  As this
holds for all $x \in Z$, we have

\begin{equation*}
\mu \Big\{ x \in X : P \big\{ \omega \in \Omega : F^*(x,\omega) > R'
(\log K)^{1/2} \big\} > 1/2 \Big\} > 4/5.
\end{equation*}

\noindent By the same Fubini trick as in the proof of Proposition 1,
we see that

\begin{gather*}
4/5 < \int_X \int_\Omega 2 \chi_{\{F^*(x,\omega) > R' (\log
K)^{1/2}\}} \, P(d\omega) \, \mu(dx) = \int_\Omega \int_X 2
\chi_{\{F^*(x,\omega) > R' (\log
K)^{1/2}\}} \, \mu(dx) \, P(d\omega) \\
= 2\left[ \int_{\mathcal{B}} \mu\{x : F^*(x,\omega) > R' (\log
K)^{1/2} \} \, P(d\omega) + \int_{\mathcal{B}^c} \mu\{x :
F^*(x,\omega) > R' (\log K)^{1/2} \} \, P(d\omega) \right] \\
\leq 2(P(\mathcal{B}) + 1/5),
\end{gather*}

\noindent which implies $P(\mathcal{B}) > 1/5$.

\vskip0.1in

\noindent {\bf Estimate of $\mathcal{C}$}.

\vskip0.1in

Now,

\begin{equation*}
\begin{split}
\int_\Omega \|\varphi(\cdot,\omega)\|_{L^1(\mu)} \, P(d\omega)
&\leq \int_\Omega \|F(\cdot,\omega)\|_{L^1(\mu)} \, P(d\omega) \\
&= \int_\Omega J^{-1/2} \int_X \left| \sum_{j=1}^J g_j(\omega) T_
j f(x) \right| \, \mu(dx) \, P(d\omega) \\
&= \int_X J^{-1/2} \int_\Omega \left| \sum_{j=1}^J g_j(\omega) T_
j f(x) \right| \, P(d\omega) \, \mu(dx) \\
&\leq \int_X J^{-1/2} \left( \int_\Omega \left| \sum_{j=1}^J
g_j(\omega)
T_j f(x) \right|^2 \, P(d\omega) \right)^{1/2} \, \mu(dx) \\
&= \int_X J^{-1/2} \left( \sum_{j=1}^J T_jf(x)^2 \right)^{1/2}
\, \mu(dx) \\
&\leq J^{-1/2} \left( \int_X \sum_{j=1}^J T_jf(x)^2 \, \mu(dx)
\right)^{1/2} \\
&= J^{-1/2} \left( \sum_{j=1}^J \|T_j f\|_{L^2(\mu)}^2 \right)^{1/2}
= J^{-1/2} \left( \sum_{j=1}^J \|f\|_{L^2(\mu)}^2 \right)^{1/2}
\\
&\leq 1.
\end{split}
\end{equation*}

\noindent From Chebyshev, we see $P\{
\|\varphi(\cdot,\omega)\|_{L^1(\mu)} > 12\} \leq 1/12 < 1/10$, or
equivalently, $P(\mathcal{C}) > 9/10$. It now follows that
$P(\mathcal{A} \cap \mathcal{B} \cap \mathcal{C}) > 0$.
\end{proof}

\section{Proofs of Theorems~\ref{thm:bj1}
and~\ref{thm:bj2}}\label{sec:banach}

Recall, our setting is a probability space $(X,\mathcal{F},\mu)$
with a sequence $S_n : L^2(\mu) \rightarrow L^2(\mu)$.
Theorem~\ref{thm:bj1} is a well-known result.  The proof is as
follows.

\begin{proof}[Proof of Theorem~\ref{thm:bj1}]
Fix $\epsilon > 0$.  Let $f \in L^p(\mu)$.  Then, $S^*f := \sup_{n}
|S_n f| < \infty$ a.s.[$\mu$], by the hypothesis.  So, there is some
$n \in \mathbb{N}$ (depending on $f$), so that $\mu \{x : S^*f(x) >
n\} \leq \epsilon/3$.  Thus,

\begin{equation*}
L^p(\mu) = \bigcup_{n = 1}^\infty \Big\{ f \in L^p(\mu) : \mu\{ S^*f
> n\} \leq \epsilon/3 \Big\}.
\end{equation*}

\noindent Denote $B_n = \{f \in L^p(\mu) : \mu\{ S^*f > n\} \leq
\epsilon/3\}$. Let $S_N^* f(x) := \sup \{ |S_n f(x)| : 1 \leq n \leq
N\}$.  It follows that

\begin{equation*}
B_n = \bigcap_{N=1}^\infty \Big\{ f \in L^p(\mu) : \mu\{ S^*_N f >
n\} \leq \epsilon/3 \Big\}.
\end{equation*}

\noindent Denote $B_n^N = \{f \in L^p(\mu) : \mu\{S^*_N f > n\} \leq
\epsilon/3\}$.  Fix $n, N$.  Now, it is clear that $\|S_N^*
f\|_{L^2(\mu)} \leq \| \sum_{k=1}^N |S_k f|\|_{L^2(\mu)} \leq
(\sum_{k=1}^N \|S_k\|) \|f\|_{L^2(\mu)} =: C\|f\|_{L^2(\mu)}$ for
all $f \in L^2(\mu)$, where $\|S_k\|$ refers to the operator norm on
$L^2(\mu)$. Suppose $(f_k) \in B_n^N$ and $f_k \rightarrow f$ in
$L^p(\mu)$-norm.  Let $r \in \mathbb{N}$. Then, by Chebyshev,

\begin{equation*}
\begin{split} \mu\{S_N^* f > n + 1/r\} &\leq \mu\{|S_N^* f - S_N^* f_k| >
1/r\} + \mu\{S_N^* f_k > n\} \\
&\leq r^2 \|S_N^*(f - f_k)\|_{L^2(\mu)}^2 + \epsilon/3 \\
&\leq C^2 r^2 \|f - f_k\|_{L^2(\mu)}^2 + \epsilon/3
\\
&\leq C^2 r^2 \|f - f_k\|_{L^p(\mu)}^2 + \epsilon/3 \rightarrow
\epsilon/3.
\end{split}
\end{equation*}

\noindent Because $\bigcup_r \{S_N^*f(x) > n + 1/r \} = \{S_N^*f(x)
> n\}$ and these sets are nested, we have the limit $\mu\{S_N^*f
> n + 1/r\} \rightarrow \mu\{S_N^*f > n\}$.  Namely, $\mu\{S_N^*f
> n\} \leq \epsilon/3$ and $f \in B_n^N$.  Thus, $B_n^N$ is
closed (in $L^p(\mu)$), which implies $B_n$ is closed.  It follows
from the Baire Category Theorem, because $L^p(\mu) = \bigcup B_n$,
that one of $B_n$ contains an open set.  That is, there exists some
$n \in \mathbb{N}$, $\delta > 0$, and $f_0 \in L^p(\mu)$ such that
$f \in B_n$ for all $\|f - f_0\|_{L^p(\mu)} \leq \delta$.  In
particular, $\mu\{ S^*(f_0 + \delta g) > n\} \leq \epsilon/3$ for
all $\|g\|_{L^p(\mu)} \leq 1$. Thus, for all $\|g\|_{L^p(\mu)} \leq
1$, we have

\begin{equation*}
\begin{split}
\mu \left\{ S^* g > \frac{2n}{\delta} \right\} \leq \mu\{ S^*(f_0 +
\delta g) > n\} + \mu\{ S^* f_0 > n\} \leq 2\epsilon/3 < \epsilon.
\end{split}
\end{equation*}

\noindent If we set $\theta(\epsilon) = 2n/\delta$, we have the
desired result.
\end{proof}

\vskip0.1in

Theorem~\ref{thm:bj2} and its proof are taken directly from Bellow
and Jones~\cite{bellowandjones}.  It is included here only for
completeness.

Denote the space $Y_0 = \{f \in L^\infty(\mu) :
\|f\|_{L^\infty(\mu)} \leq 1\}$.  We will be concerned with the
$L^2(\mu)$-norm on $Y_0$.  It is well-known that $Y_0$ is complete
under $\| \cdot \|_{L^2(\mu)}$.  Denote $\delta$-balls in $Y_0$ by
$B_\delta(f) = \{g \in Y_0 : \|f - g\|_{L^2(\mu)} < \delta\}$. The
first step is to prove the following lemma.

\begin{lemma}\label{lemma:bj3}
If $f_0 \in Y_0$ and $\delta > 0$, then $B_\delta(0) \subseteq
B_\delta(f_0) - B_\delta(f_0) = \{g_1 - g_2 : g_1, g_2 \in
B_\delta(f_0)\}$.  \end{lemma}

\begin{proof} Let $g \in B_\delta(0)$, so that $\|g\|_{L^2(\mu)} <
\delta$.  Define

\begin{equation*}
\begin{split}
u_1(x) &= \begin{cases} g(x)& \quad \text{ if } f_0(x), g(x) \,\,
\text{are both finite, and } f_0(x)g(x) \leq 0,
\\
0& \quad \text{ otherwise,} \end{cases} \\
u_2(x) &= \begin{cases} -g(x)& \,\, \text{ if } f_0(x), g(x) \,\,
\text{are both finite, and } f_0(x)g(x) > 0,
\\
0& \,\, \text{ otherwise.} \end{cases}
\end{split}
\end{equation*}

\noindent Let $g_1 = f_0 + u_1$ and $g_2 = f_0 + u_2$.  Now, it is
clear that $g = u_1 - u_2 = (f_0 + u_1) - (f_0 + u_2) = g_1 - g_2$
for a.s.[$\mu$] $x \in X$.  It is also clear that $\|f_0 -
g_1\|_{L^2(\mu)} = \|u_1\|_{L^2(\mu)} \leq \|g\|_{L^2(\mu)} <
\delta$.  Similarly, $\|f_0 - g_2\|_{L^2(\mu)} < \delta$.  So, we
will be done if we can show $g_1, g_2 \in Y_0$.

Fix an $x$ such that $f_0(x), g(x)$ are both finite and $f_0(x)g(x)
\leq 0$, i.e., $f_0(x)$ and $g(x)$ have opposite signs.  Then,
$u_1(x) = g(x)$ and $|g_1(x)| = |f_0(x) + u_1(x)| = |f_0(x) +
g(x)|$.  As they have opposite signs, $|f_0(x) + g(x)| \leq \max\{
|f_0(x)|, |g(x)|\}$. Now suppose $x$ is such that $f_0(x), g(x)$ are
finite and $f_0(x)g(x) > 0$. Then, $u_1(x) = 0$ and $|g_1(x)| =
|f_0(x)|$.  Hence, $|g_1| \leq \max\{ |f_0|, |g| \}$ for a.s.[$\mu$]
$x$, which implies $\|g_1\|_{L^\infty(\mu)} \leq \max\{
\|f_0\|_{L^\infty(\mu)}, \|g\|_{L^\infty(\mu)}\} \leq 1$.  Namely,
$g_1 \in Y_0$.  Precisely the same argument shows $g_2 \in Y_0$.
\end{proof}

\vskip0.1in

We can now proceed to the proof of Theorem~\ref{thm:bj2}.  This
proof also relies on the Baire Category Theorem.

\begin{proof}[Proof of Theorem~\ref{thm:bj2}]
Fix $\epsilon , \eta > 0$.  Choose $0 < \alpha < 1/2$ so that
$\alpha < \epsilon/3$ and $\eta > 2\alpha$.  For $N \in \mathbb{N}$
define

\begin{equation*}
F_N(\alpha) = \left\{ f \in Y_0 : \mu \left\{x \in X : \sup_{m \geq
N} |S_N f(x) - S_m f(x)| \leq \alpha \right\} \geq 1 - \alpha
\right\}.
\end{equation*}

\noindent For $M > N$, define

\begin{equation*}
F_{N,M}(\alpha) = \left\{ f \in Y_0 : \mu \left\{x \in X : \sup_{N
\leq m \leq M} |S_N f(x) - S_m f(x)| \leq \alpha \right\} \geq 1 -
\alpha \right\}.
\end{equation*}

Fix $N$ and $M > N$.  We wish to show $F_{N,M}(\alpha)$ is closed,
with respect to $(Y_0, \|\cdot\|_{L^2(\mu)})$. Let $(f_k) \in
F_{N,M}(\alpha)$ and $f_k \rightarrow f \in Y_0$ in $L^2(\mu)$-norm.
Let $g_k(x) = \sup_{N \leq m \leq M} |S_N f_k(x) - S_m f_k(x)|$ and
$g(x) = \sup_{N \leq m \leq M} |S_N f(x) - S_m f(x)|$.  By $\|S_n\|$
it is meant the operator norm on $L^2(\mu)$.  Now,

\begin{equation*}
\begin{split}
\|g - g_k\|_{L^2(\mu)} &\leq \left\| \sup_{N \leq m \leq M} |S_N f -
S_m f - S_N f_k + S_m f_k| \right\|_{L^2(\mu)} \\
&\leq \|S_N(f - f_k)\|_{L^2(\mu)} + \left\| \sup_{N \leq m \leq M}
|S_m(f_k - f)| \right\|_{L^2(\mu)} \\
&\leq \|S_N(f - f_k)\|_{L^2(\mu)} + \left\| \sum_{m = N}^M |S_m(f_k
- f)| \right\|_{L^2(\mu)} \\
&\leq \left(\|S_N\| + \sum_{k=N}^M \|S_k\| \right) \|f -
f_k\|_{L^2(\mu)} \rightarrow 0.
\end{split}
\end{equation*}

\noindent For $n \in \mathbb{N}$, we have that

\begin{equation*}
\begin{split}
\mu \left\{x \in X : g(x) > \alpha + \frac{1}{n} \right\} &\leq \mu
\left\{x \in X : g_k(x) > \alpha \right\} + \mu \left\{x \in X : |g
- g_k| > \frac{1}{n} \right\} \\
&\leq \alpha + n^2 \|g - g_k\|_{L^2(\mu)}^2 \rightarrow \alpha.
\end{split}
\end{equation*}

\noindent As $\{g > \alpha\} = \bigcup_n \{g > \alpha +
\frac{1}{n}\}$, and this is a nested sequence, we see $\mu\{g >
\alpha\} \leq \alpha$.  But this says $f \in F_{N,M}(\alpha)$, and
$F_{N,M}(\alpha)$ is closed.

Now, as the relevant sets are nested and decreasing, it is easy to
see that $F_N(\alpha) = \bigcap_{M > N} F_{N,M}(\alpha)$.  So, each
$F_N(\alpha)$ is closed.  Let $g \in Y_0$.  By hypothesis, $S_n g$
converges a.s.[$\mu$].  By Egoroff's Theorem, there is a set $E$
with $\mu(X - E) < \alpha$ so that $S_n g$ converges uniformly on
$E$.  Then, there is some $N$ such that $|S_N g(x) - S_n g(x)| \leq
\alpha$ whenever $n \geq N$ and $x \in E$.  This implies $g \in
F_N(\alpha)$.  Hence, $Y_0 = \bigcup_N F_N(\alpha)$.  As $Y_0$ is
complete under $\|\cdot\|_{L^2(\mu)}$, it follows by the Baire
Category Theorem that at least one of $F_N(\alpha)$ contains an open
set. That is, there is some $N_0$, some $f_0 \in Y_0$, and some
$\delta > 0$ such that $B_\delta(f_0) \subset F_{N_0}(\alpha)$.

Let $S^*_{N_0} f = \sup_{1 \leq n \leq N_0} |S_n f|$ and $S^* f =
\sup_n |S_n f|$.  Note, for each $g \in Y_0$, $\|S_n g\|_{L^2(\mu)}
\leq \|S_n\| \|g\|_{L^2(\mu)}$, so that each $S_n$ is continuous at
0 in $(Y_0, \|\cdot\|_{L^2(\mu)})$.  Hence, $S^*_{N_0}$ is
continuous at 0, because $S_{N_0}^*(f) \leq \sum_{n=1}^{N_0} |S_n
f|$. So, there is some $0 < \delta' < \delta$ such that $\|S^*_{N_0}
f\|_{L^2(\mu)}^2 < \alpha(\eta - 2\alpha)^2$ when $f \in
B_{\delta'}(0)$.

Fix $f \in B_{\delta'}(0)$.  Now, by Lemma~\ref{lemma:bj3}, we see
$B_{\delta'}(0) \subset B_{\delta}(0) \subseteq B_{\delta}(f_0) -
B_{\delta}(f_0)$. Thus, there are $g_1, g_2 \in B_\delta(f_0)
\subset F_{N_0}(\alpha)$ such that $f = g_1 - g_2$ a.s.[$\mu$]. By
definition of $F_{N_0}(\alpha)$, we see that

\begin{gather*}
\mu \left\{ x \in X : \sup_{m \geq N_0} |S_{N_0} g_1(x) - S_m
g_1(x)| \leq \alpha \right\} \geq 1 - \alpha, \\
\mu \left\{ x \in X : \sup_{m \geq N_0} |S_{N_0} g_2(x) - S_m
g_2(x)| \leq \alpha \right\} \geq 1 - \alpha.
\end{gather*}

\noindent Now, except for a set of probability 0,

\begin{equation*}
\begin{split}
\sup_{m \geq N_0} \big| S_{N_0} f(x) - S_m f(x) \big| &= \sup_{m
\geq N_0} \big| S_{N_0} g_1(x) - S_m g_1(x) - S_{N_0} g_2(x) +
S_m g_2(x) \big| \\
&\leq \sup_{m \geq N_0} \big| S_{N_0} g_1(x) - S_m g_1(x) \big| +
\sup_{m \geq N_0} \big| S_{N_0} g_2(x) - S_m g_2(x) \big|.
\end{split}
\end{equation*}

\noindent Therefore,

\begin{equation*}
\begin{split}
\mu \bigg\{ x \in X : &\sup_{m \geq N_0} |S_{N_0} f(x) - S_m f(x)|
\leq 2\alpha \bigg\} \\
&\geq \mu \bigg( \Big\{ \sup_{m \geq N_0} |S_{N_0} g_1 - S_m g_1|
\leq \alpha \Big\} \bigcap \Big\{ \sup_{m
\geq N_0} |S_{N_0} g_2 - S_m g_2| \leq \alpha \Big\} \bigg) \\
&\geq 1 - 2 \alpha.
\end{split}
\end{equation*}

Define

\begin{gather*}
C = \Big\{x \in X : S^* f(x) > \eta \Big\}, \\
D = \Big\{x \in X : S^*f(x) \leq S^*_{N_0}f(x) + 2\alpha \Big\}, \\
E = \Big\{x \in X: \sup_{m \geq N_0} |S_{N_0} f(x) - S_m f(x)| \leq
2\alpha \Big\}.
\end{gather*}

\noindent Then, $\mu(E) \geq 1 - 2\alpha$ by above. Now, for any $x
\in X$ and $m \leq N_0$ it is clear that

\begin{equation*}
|S_m f(x)| \leq S^*_{N_0} f(x).
\end{equation*}

\noindent On the other hand, for $x \in E$ and any $m \geq N_0$ we
have

\begin{equation*}
|S_m f(x)| \leq |S_{N_0} f(x)| + 2\alpha \leq S^*_{N_0} f(x) +
2\alpha.
\end{equation*}

\noindent Namely, for every $x \in E$ we see $S^* f(x) \leq
S^*_{N_0} f(x) + 2\alpha$, or $x \in D$.  Thus, $\mu(D) \geq \mu(E)
\geq 1 - 2\alpha$.  Finally,

\begin{equation*}
\begin{split}
\mu(C) &= \mu(C \cap D) + \mu(C \cap D^c) \\
&\leq \mu \{x \in X : S^*_{N_0} f(x) > \eta - 2\alpha\} + \mu(D^c)
\\
&\leq \|S^*_{N_0} f\|_{L^2(\mu)}^2 \frac{1}{(\eta - 2\alpha)^2} +
2\alpha
\\
&< 3\alpha < \epsilon.
\end{split}
\end{equation*}

\noindent This holds for all $f \in B_{\delta'}(0)$.  If we set
$\rho(\epsilon, \eta) = \delta'$, we have the desired result.
\end{proof}

\section{Proof of Theorem~\ref{thm:fernique}}\label{sec:fern}

For this section, it will be convenient to discuss normal random
vectors, that is, measurable maps $G : \Omega \rightarrow
\mathbb{R}^N$ where $G = (G_1, \ldots, G_N)$ and each $G_j$ is a
normal random variable.  We will say $G$ has mean 0 if each $G_j$
has mean 0.

The concepts of distribution and independence are easily extended to
this case.  We say two random vectors $G$ and $H$ have the same
distribution if $P(G \in E) = P(H \in E)$ for all measurable sets $E
\subset \mathbb{R}^N$. We say $G$ and $H$ are independent if $P(G
\in E) P(H \in F) = P(G \in E, H \in F)$ for all $E, F \subset
\mathbb{R}^N$.

\begin{lemma}\label{lemma:fern}
Let $G$ and $H$ be mean 0 normal random vectors which are
independent and have the same distribution. Then, for any $\theta
\in \mathbb{R}$, the normal random vector $(G,H) : \Omega
\rightarrow \mathbb{R}^{2N}$ has the same distribution as $(G \sin
\theta + H \cos \theta, G \cos \theta - H \sin \theta)$.
\end{lemma}

\begin{proof} A normal random vector (and its distribution) is
completely determined by its covariance matrix $\cov(s,t) = E(G_s
G_t)$.  But, it is clear that for all $s, t$

\begin{gather*}
E(G_s G_t) = E((G_s \sin \theta + H_s \cos \theta)(G_t
\sin \theta + H_t \cos \theta)), \\
E(H_s H_t) = E((G_s \cos \theta - H_s \sin \theta)(G_t
\cos \theta - H_t \sin \theta)), \\
E(G_s H_t) = 0 = E((G_s \sin \theta + H_s \cos \theta)(G_t \cos
\theta - H_t \sin \theta)).
\end{gather*}

\noindent Hence, the covariance matrices of $(G,H)$ and $(G \sin
\theta + H \cos \theta, G \cos \theta - H \sin \theta)$ are the
same.
\end{proof}

\vskip0.1in

The above lemma is, in some sense, the statement that mean 0 normal
random vectors are rotation invariant.

Now, for a mean 0 normal random vector $G$ on $(\Omega, \mathcal{B},
P)$, we can always find normal random vectors $H$ and $K$ on some
probability space $(\Omega', \mathcal{B}', P')$ such that $H, K$ are
independent, and $H, K$ have the same distribution as $G$.  That is
$P(G \in E) = P'(H \in E) = P'(K \in E)$ for all measurable sets $E
\subset \mathbb{R}^N$.  We say $H, K$ are independent copies of $G$.
We now prove Theorem~\ref{thm:fernique}.  This proof is taken
from~\cite{fernique}.

\begin{proof} [Proof of Theorem~\ref{thm:fernique}]
Write $G = (G_1, \ldots, G_N)$ as a normal random vector.  Let $H,
K$ be independent copies of $G$ on some probability space $(\Omega',
\mathcal{B}', P')$. Define $S(G) : \Omega \rightarrow [0,\infty)$ by
$S(G)(\omega) = \sup_n |G_n(\omega)|$. Define $S(H), S(K)$
similarly.  It follows $S(G), S(H), S(K)$ have the same distribution
and that $S(H), S(K)$ are independent. Then, for $t \geq s$

\begin{equation*}
\begin{split}
2 P\{S(G) \leq s\} P\{S(G) > t\} &= P'\{S(H) \leq s\} P'\{S(K) > t\}
+
P'\{S(K) \leq s\} P'\{S(H) > t\} \\
&= P'\{S(H) \leq s, S(K) > t\} + P'\{S(K) \leq s, S(H) > t\}. \\
\end{split}
\end{equation*}

\noindent Apply Lemma~\ref{lemma:fern} to $H, K$ with $\theta =
\pi/4$. Then, $(H,K)$ and $(\frac{H+K}{\sqrt{2}},
\frac{H-K}{\sqrt{2}})$ have the same distribution.  It follows that
$P'\{S(H) \leq s, S(K) > t\} = P'\{S(\frac{H+K}{\sqrt{2}}) \leq s,
S(\frac{H-K}{\sqrt{2}}) > t\}$. Similarly for the second term, so
that

\begin{equation*}
\begin{split}
&2P\{S(G) \leq s\} P\{S(G) > t\} \\
&= P' \left\{ S\left(\frac{H+K}{\sqrt{2}} \right) \leq s,
S\left(\frac{H-K}{\sqrt{2}}\right)
> t \right\} + P'\left\{S \left(\frac{H-K}{\sqrt{2}} \right) \leq
s, S\left(\frac{H+K}{\sqrt{2}}\right) > t \right\} \\
&= P' \left( \Big\{ S(H+K) \leq s\sqrt{2}, S(H-K)
> t\sqrt{2} \Big\} \bigcup \Big\{ S(H-K) \leq
s\sqrt{2}, S(H+K) > t\sqrt{2} \Big\} \right),
\end{split}
\end{equation*}

\noindent where the last equality follows as the sets are clearly
disjoint.  Consider the first set in the union, $\{S(H+K) \leq
s\sqrt{2}, S(H-K) > t\sqrt{2} \}$.  In this set $\sqrt{2}(t-s) <
S(H+K) - S(H-K) \leq S(2K) = 2S(K)$ or $\sqrt{2}S(K) > t - s$.
Similarly, $\sqrt{2}(t-s) < S(H+K) - S(K-H) \leq S(2H) = 2S(H)$ or
$\sqrt{2}S(H) > t-s$.  The same calculations work in the other set.
So, we have

\begin{equation}\label{eq:fern1}
\begin{split}
2P\{S(G) \leq s\} P\{S(G) > t\} &\leq P'\{\sqrt{2}S(H) > t-s,
\sqrt{2}S(K) > t-s\} \\
&= P\{\sqrt{2}S(G) > t-s\}^2.
\end{split}
\end{equation}

We now define a sequence $(t_n)$.  Set $t_0 = s$ and $t_{n+1} =
t_n\sqrt{2} + s$.  It is easily checked by induction that $t_n =
(\sqrt{2} + 1)(2^{(n+1)/2} - 1)s$.  Define $q = 2P\{S(G) \leq s\}$
and $x_n = q^{-1} P\{S(G) > t_n\}$.  By the hypothesis, $q \geq 1$.
By construction and from (\ref{eq:fern1}),

\begin{equation*}
\begin{split}
q^2 x_{n+1} &= q P\{S(G) > t_{n+1}\} = 2 P\{S(G) \leq s\} P\{S(G) >
t_{n+1}\} \\
&\leq P \left\{S(G) > \frac{t_{n+1} - s}{\sqrt{2}} \right\}^2 =
P\{S(G)
> t_n\}^2 = q^2 x_n^2,
\end{split}
\end{equation*}

\noindent or $x_{n+1} \leq x_n^2$.  This implies that $x_n \leq
x_0^{2^n}$. It is easily seen that $x_0 = \frac{2-q}{2q} \leq 1/2$
as $q \geq 1$. Hence, $P\{S(G) > t_n\} = qx_n \leq q x_0^{2^n} \leq
q 2^{-2^n}$. By Lemma~\ref{lemma:expect},

\begin{equation*}
\begin{split}
\|S(G)\|_{L^1(P)} &= \int_0^\infty P\{S(G) > t\} \, dt \\
&= \int_0^s P\{S(G) > t\} \, dt + \sum_{n=0}^\infty
\int_{t_n}^{t_{n+1}} P\{S(G) > t\} \, dt \\
&\leq s + \sum_{n=0}^\infty \int_{t_n}^{t_{n+1}} P\{S(G) > t_n\} \,
dt \\
&\leq s + \sum_{n=0}^\infty (t_{n+1} - t_n) q2^{-2^n}.
\end{split}
\end{equation*}

\noindent Of course, $q = 2P\{S(G) \leq s\} \leq 2$ trivially.  From
the induction characterization of $t_n$, we see $t_{n+1} - t_n =
s2^{(n+1)/2}$.  Therefore,

\begin{equation*}
\begin{split}
\|S(G)\|_{L^1(P)} &\leq s + s2^{3/2} \sum_{n=0}^\infty 2^{n/2}
2^{-2^n}
\leq s + s 2^{3/2} \sum_{n=0}^\infty 2^{n/2} 2^{-2n} \\
&= s + s 2^{3/2} \sum_{n=0}^\infty 2^{-3n/2} = s + s
\frac{2^{3/2}}{1- 2^{-3/2}} \leq 6s.
\end{split}
\end{equation*}
\end{proof}

\section{An Application}

Let $\mathbf{T} = \mathbb{R}/\mathbb{Z}$.  A function $f$ on
$\mathbb{R}$ with period 1 can be viewed as a function on
$\mathbf{T}$.  Let $m$ be Lebesgue measure, and consider the
probability space $(\mathbf{T},m)$.  Let $(a_j)$ be any non-zero
sequence of real numbers which converge to 0. For $f : \mathbf{T}
\rightarrow \mathbb{R}$, consider the operators

\begin{equation*}
S_n f(x) = \frac{1}{n} \sum_{j=1}^n f(x+a_j).
\end{equation*}

\noindent Bellow asked whether $S_n f$ converges to $f$ a.s.~for all
$f \in L^1(m)$.  The answer to this turns out to be no. In fact, it
is not even true for all $f \in L^\infty(m)$.

We will prove this using the second entropy result. In this case, it
is beneficial to consider complex-valued functions temporarily.
Here, the operators $S_n$ make perfectly good sense applied to
complex-valued functions.  In fact, we also have the nice property
that $S_n(\re f) = \re(S_n f)$ and similarly for the imaginary part.
We will take advantage of this.  First, we need two technical
lemmas.

\begin{lemma}\label{lemma:tech} Let $(a_j)$ be a sequence of
non-zero real numbers converging to 0.  Then, given any $r \in
\mathbb{N}$, there exist integers $J_1 < J_2 < \ldots < J_r$
satisfying the following: if $\bar{\alpha} = (\alpha_1, \ldots,
\alpha_r)$ is a vector of 0's and 1's, then there is an integer
$n(\bar{\alpha})$ such that

\begin{gather*}
\left| 1 - J_s^{-1} \sum_{j \leq J_s} e^{2\pi ia_j n(\bar{\alpha})}
\right| < \frac{1}{10} \quad \text{if } \,\, \alpha_s = 0, \\
\left| 1 - J_s^{-1} \sum_{j \leq J_s} e^{2\pi ia_j n(\bar{\alpha})}
\right| > \frac{1}{2} \quad \text{if } \,\, \alpha_s = 1, \\
\end{gather*}

\noindent for all $1 \leq s \leq r$.
\end{lemma}

\begin{lemma}\label{lemma:comb}
Let $(Y, \mathcal{S},\nu)$ be a measure space and define $\|h\| =
(\int_Y |h|^p \, d\nu)^{1/p}$ for the measurable functions from $Y$
to $\mathbb{C}$, where $1 \leq p < \infty$. Let $r_0 \in \mathbb{N}$
and $r = 4r_0^2 + 2r_0 \in \mathbb{N}$. Suppose $\{h_1, \ldots,
h_r\}$ is a collection of complex-valued functions on $Y$ such that
$\|h_j - h_k\| > \alpha$ for all $j \not= k$. Then, there exists a
set $I \subset \{1, \ldots, r\}$, $|I| = r_0$ such that either
$\|\re h_j - \re h_k\| > \alpha/4$ for all $j \not= k \in I$ or
$\|\im h_j - \im h_k\| > \alpha/4$ for all $j \not= k \in I$.
\end{lemma}

We temporarily postpone the proofs of these lemmas and proceed to
the solution of Bellow's question.

\begin{thm} Let $(a_j)$ be any real sequence of numbers which
converge to 0 and $a_j \not= 0$ for all $j$.  Then, there exists $f
\in L^\infty(m)$ such that $S_n f$ does not converge a.s.[$m$].
\end{thm}

\begin{proof}
Let $T_j f(x) = f(x + b_j)$ for some real sequence $(b_j)$.  Then,
the fact that $T_j(1) = 1$ and $T_j$ are positive is obvious. By
periodicity, $T_j$ is an isometry on $L^1(m)$ and $L^2(m)$. Also,
$\|S_n f\|_{L^2(m)} \leq \frac{1}{n} \sum \|f\|_{L^2(m)} =
\|f\|_{L^2(m)}$ and $T_j S_n = S_n T_j$.

Let $w$ be an irrational number and $b_j = (j-1)w$. It then follows
from the equidistribution theorem (or a special case of Birkhoff's
Ergodic Theorem) that $(T_j)$ satisfies the mean ergodic condition.
Thus, $(S_n)$ commutes with a Bourgain sequence, and the second
entropy result can be applied to $(S_n)$.

It suffices to show that for some $\delta > 0$ we have $\sup
\{N_f(\delta) : \|f\|_{L^2(m)} \leq 1\} = \infty$.  In fact, we will
do this with $\delta = 1/40$.  Let $r_0 \in \mathbb{N}$ and $r =
4r_0^2 + 2r_0$. By Lemma~\ref{lemma:tech}, choose integers $J_1 <
\ldots < J_r$. Define a complex-valued function $g : \mathbf{T}
\rightarrow \mathbb{C}$ by

\begin{equation*}
g(x) = 2^{-r/2} \sum_{\bar{\alpha} \in \{0,1\}^r} e^{2\pi i
n(\bar{\alpha})x}.
\end{equation*}

\noindent  As we said before, we can consider $S_n$ acting on
complex-valued functions.  Although the second entropy result cannot
be applied in this case, we will use $g$ to manufacture an
appropriate real-valued function.  Note, $\|g\|_2^2 = \int_0^1 g(x)
\bar{g}(x) \, dx = 2^{-r} \sum_{\bar{\alpha}} 1 = 1$ by
orthogonality.  (The notation $\|\cdot\|_2$ has the obvious meaning,
where we make the distinction here from $L^2(m)$ and
$\|\cdot\|_{L^2(m)}$ which implies real-valued functions).  Also,

\begin{equation*}
S_{J_s} g(x) = 2^{-r/2} \sum_{\bar{\alpha} \in \{0,1\}^r}
\beta_{s,\bar{\alpha}} e^{2\pi i n(\bar{\alpha}) x},
\end{equation*}

\noindent where

\begin{equation*}
\beta_{s,\bar{\alpha}} = J_s^{-1} \sum_{j \leq J_s} e^{2\pi ia_j
n(\bar{\alpha})}.
\end{equation*}

\noindent Fix an $\bar{\alpha}$ and suppose $\alpha_s = 1$ and
$\alpha_t = 0$.  By Lemma~\ref{lemma:tech}, we have

\begin{equation*}
|\beta_{s,\bar{\alpha}} - \beta_{t,\bar{\alpha}}| \geq
|\beta_{s,\bar{\alpha}} - 1| - |1 - \beta_{t,\bar{\alpha}}| >
\frac{1}{2} - \frac{1}{10} = \frac{2}{5}.
\end{equation*}

\noindent By symmetry, this holds so long as $\alpha_s \not=
\alpha_t$.  Hence, by orthogonality and above,

\begin{equation*}
\begin{split}
\|S_{J_s} g - S_{J_t} g\|_2 &= 2^{-r/2} \left( \sum_{\bar{\alpha}}
|\beta_{s,\bar{\alpha}} - \beta_{t,\bar{\alpha}}|^2 \right)^{1/2}
\\
&\geq 2^{-r/2} \left( \sum_{\bar{\alpha} : \alpha_s \not= \alpha_t}
|\beta_{s,\bar{\alpha}} - \beta_{t,\bar{\alpha}}|^2 \right)^{1/2}
\\
&> 2^{-r/2} (2/5) \left( \sum_{\bar{\alpha} : \alpha_s \not=
\alpha_t} 1 \right)^{1/2} \\
&= 2^{-r/2} (2/5) 2^{(r-1)/2} \\
&> 1/5.
\end{split}
\end{equation*}

\noindent This holds for all $s \not= t$.

Apply Lemma~\ref{lemma:comb} to the set $\{S_{J_1} g, \ldots,
S_{J_r} g\}$ to find a subset $I \subset \{J_1, \ldots, J_r\}$, $|I|
= r_0$ such that either $\|\re S_{J_t} g - \re S_{J_s} g\|_2 > 1/20$
or $\|\im S_{J_t} g - \im S_{J_s} g\|_2 > 1/20$ for all $J_s \not=
J_t \in I$. If it is the first, set $f = \re g$, and if it is the
second, set $f = \im g$. Then, $\|f\|_{L^2(m)} \leq \|g\|_2 = 1$.
Further, $\|S_{J_s} f - S_{J_t} f\|_{L^2(m)} > 1/20$ for all $J_s
\not= J_t \in I$.  As no two such $S_{J_s} f$ could be contained in
the same $1/40$-ball in $L^2(m)$, we see $N_f(1/40) \geq |I| = r_0$.
As $r_0$ is arbitrary, $\sup N_f(1/40) = \infty$.
\end{proof}

To conclude, we need only establish Lemmas~\ref{lemma:tech}
and~\ref{lemma:comb}. First, we recall three simple results in
complex arithmetic.

\begin{claim} Let $a, b \in \mathbb{C}$, with $|a|, |b| \leq 1$, and
$\lambda \in \mathbb{R}$.  Then,

\begin{enumerate}
\item $|ab - 1| \leq |a - 1| + |b - 1|$,  \item $\re (ab) \leq |a-1|
+ \re(b)$, \item $|1 - e^{2\pi i\lambda}| \leq 2\pi|\lambda|$.
\end{enumerate}
\end{claim}

\begin{proof} First, $|ab - 1| = |ab - a + a - 1| \leq |ab - a| +
|a - 1| \leq |b-1| + |a-1|$.  Second, $\re(ab) = \re(ab - b) +
\re(b) \leq |ab - b| + \re(b) \leq |a - 1| + \re(b)$.  Third, recall
$1 - \cos(2x) \leq 2x^2$ for all real $x$. So, $|1 - e^{2\pi
i\lambda}|^2 = (1 - \cos(2\pi\lambda))^2 + \sin^2(2\pi\lambda) = 2(1
- \cos(2\pi \lambda)) \leq 4\pi^2 \lambda^2$.
\end{proof}

\vskip0.2in

\begin{proof}[Proof of Lemma~\ref{lemma:tech}]
Fix $r \in \mathbb{N}$.  If $r = 1$, then set $J_1 = 1$ and choose
$n(\bar{\alpha})$ accordingly. Assume $r > 1$.  For each $1 \leq s
\leq r$, we will construct integers $m_s$ simultaneously as $J_s$.

Set $J_1 = 1$ and choose $m_1$ so that $|1 - e^{2\pi ia_1 m_1}| >
3/4$.  Assume $J_t$ and $m_t$ are known for all $t < s$.  Let $M_s =
\sum_{t < s} |m_t|$. As $a_j \rightarrow 0$, there is a $L_s > 0$
such that $\sup_{j > L_s/100} |a_j| \leq (400M_s\pi)^{-1}$. Further,
we can choose $T_s > 0$ (depending on $a_j$ for $j \leq J_{s-1}$)
such that for each $z \in \mathbb{Z}$ there is a corresponding $t
\in \mathbb{Z}$, $|t| \leq T_s$ satisfying $|e^{2\pi ia_jz} -
e^{2\pi ia_jt}| < 1/50r$ for all $j \leq J_{s-1}$.

As $a_j \rightarrow 0$, we can choose $J_s$ such that $J_s > L_s$,
$J_s > J_{s-1}$, and $J_s^{-1} \sum_{j \leq J_s} |a_j| <
(100T_s)^{-1}$. Also, as $a_j \not= 0$ for all $j$, note that

\begin{equation*}
\begin{split}
\lim_{R \rightarrow \infty} \left| \frac{1}{R} \int_0^R \left(
\frac{1}{J_s} \sum_{j \leq J_s} e^{2\pi i a_jx} \right) \, dx
\right| &= \lim_{R \rightarrow \infty} \left| \frac{1}{R}
\left(\frac{1}{J_s}
\sum_{j \leq J_s} \frac{1}{2\pi i a_j} (e^{2\pi ia_j R}  - 1) \right) \right| \\
&\leq \lim_{R \rightarrow \infty} \frac{1}{R} J_s^{-1} \sum_{j \leq
J_s} \frac{1}{\pi a_j} = 0.
\end{split}
\end{equation*}

\noindent It follows there is $y_s > 0$ such that $\re(J_s^{-1}
\sum_{j \leq J_s} e^{2\pi ia_j y_s}) < 1/10$; otherwise, this limit
could not be 0. Set $z_s$ to be the integer part of $y_s$, and take
$|t_s| \leq T_s$ as prescribed above. Set $m_s = z_s - t_s$. Define
all $J_s$ and $m_s$ in this manner.

Then, for each $1 < s \leq r$, by second and third statements of the
above claim and by construction, we have

\begin{equation*}
\begin{split}
\re\left(J_s^{-1} \sum_{j \leq J_s} e^{2\pi ia_jm_s} \right) &=
J_s^{-1} \sum_{j \leq J_s} \re(e^{2\pi ia_jy_s} e^{2\pi ia_j(m_s -
y_s)})
\\
&\leq J_s^{-1} \sum_{j \leq J_s} \left (\re(e^{2\pi ia_jy_s}) + |1 -
e^{2\pi ia_j(m_s - y_s)}| \right) \\
&= \re \left(J_s^{-1} \sum_{j \leq J_s} e^{2\pi ia_j y_s} \right) +
J_s^{-1} \sum_{j \leq J_s} |1 - e^{2\pi ia_j(z_s - y_s - t_s)}| \\
&\leq \frac{1}{10} + 2\pi J_s^{-1} \sum_{j \leq J_s} |a_j|(|z_s -
y_s| + |t_s|) \\
&\leq \frac{1}{10} + \frac{2\pi (T_s+1)}{100 T_s} \leq \frac{1}{10}
+ \frac{4\pi}{100} < 1/4.
\end{split}
\end{equation*}

\noindent This gives

\begin{gather*}
J_s > L_s, \\
|1 - e^{2\pi ia_jm_s}| = |e^{2\pi ia_jz_s} - e^{2\pi ia_jt_s}| <
1/50r
\quad \text{ for all } j \leq J_{s-1}, \\
\left| 1 - J_s^{-1} \sum_{j \leq J_s} e^{2\pi ia_j m_s} \right| \geq
1 - \re \left( J_s^{-1} \sum_{j \leq J_s} e^{2\pi ia_j m_s} \right)
> 3/4.
\end{gather*}

\noindent But, $J_1 = 1$ and $m_1$ was chosen so the last condition
is also true for $s = 1$.  Further, if we rewrite the second
condition, we have

\begin{gather*}
\sup_{j \leq J_{t-1}} |1 - e^{2\pi ia_j m_t}| < 1/50r \quad \text{
for all } 2 \leq t \leq r, \\
\left| 1 - J_s^{-1} \sum_{j \leq J_s} e^{2\pi ia_j m_s} \right|
> 3/4 \quad \text{ for all } 1 \leq s \leq r.
\end{gather*}

Fix $\bar{\alpha} \in \{0,1\}^r$.  Define $n_s = \alpha_s m_s$ and
$n(\bar{\alpha}) = n_1 + \ldots + n_r$. Fix $1 \leq s \leq r$. Then,
by the first statement in the claim,

\begin{gather*}
\left| J_s^{-1} \sum_{j \leq J_s} e^{2\pi ia_j n(\bar{\alpha})} -
J_s^{-1} \sum_{j \leq J_s} e^{2\pi i a_jn_s} \right| \leq J_s^{-1}
\sum_{j \leq J_s} | e^{2\pi ia_jn(\bar{\alpha})} - e^{2\pi ia_jn_s} | \\
= J_s^{-1} \sum_{j \leq J_s} |e^{2\pi ia_j n_s}| |e^{2\pi
ia_j(n(\bar{\alpha}) - n_s)} - 1| \\
= J_s^{-1} \sum_{j \leq J_s} \left| \exp \bigg(2\pi ia_j
\Big(\sum_{t < s} n_t \Big) \bigg) \exp \bigg(2\pi ia_j
\Big(\sum_{t > s} n_t \Big) \bigg) - 1 \right| \\
\leq J_s^{-1} \sum_{j \leq J_s} \left| \exp \bigg(2\pi ia_j
\Big(\sum_{t < s} n_t \Big) \bigg) - 1 \right| + J_s^{-1} \sum_{j
\leq J_s} \left| \exp \bigg(2\pi ia_j \Big(\sum_{t > s} n_t \Big)
\bigg) - 1 \right| =: I + II.
\end{gather*}

\noindent If $s = 1$, then $I = 0$.  If $s > 1$, then as $J_s >
L_s$,

\begin{equation*}
\begin{split}
I &= J_s^{-1} \sum_{j < J_s/100} \left| \exp \bigg(2\pi ia_j
\Big(\sum_{t < s} n_t \Big) \bigg) - 1 \right| + J_s^{-1}
\sum_{J_s/100 \leq j \leq J_s} \left| \exp \bigg(2\pi ia_j
\Big(\sum_{t < s} n_t \Big) \bigg) - 1 \right| \\
&\leq J_s^{-1} \left( \sum_{j < J_s/100} 2 \right) + J_s^{-1} \left(
\sum_{J_s/100 \leq j \leq J_s} 2\pi |a_j| \Big| \sum_{t <
s} n_s \Big| \right) \\
&\leq \frac{1}{50} + 2M_s\pi J_s^{-1} (J_s - J_s/100) \sup_{j >
J_s/100}
|a_j| \\
&\leq \frac{1}{50} + 2M_s\pi \sup_{j > L_s/100} |a_j| < \frac{1}{50}
+ \frac{1}{200} = \frac{1}{40}.
\end{split}
\end{equation*}

\noindent On the other hand, for all $s$, we have by applying the
first statement of the claim several times

\begin{equation*}
\begin{split}
II &= J_s^{-1} \sum_{j \leq J_s} \left| \exp \bigg(2\pi ia_j
\Big(\sum_{t > s} n_t \Big) \bigg) - 1 \right| \\
&\leq \sup_{j \leq J_s} \left| \exp \bigg(2\pi ia_j \Big(\sum_{t >
s} n_t \Big) \bigg) - 1 \right| \\
&= \sup_{j \leq J_s} \left| \prod_{t > s} e^{2\pi ia_jn_t} - 1
\right| \leq
\sup_{j \leq J_s} \sum_{t > s} |e^{2\pi ia_jn_t} - 1| \\
&\leq \sum_{t > s} \sup_{j \leq J_s} |e^{2\pi ia_jm_t} - 1| \leq
\sum_{t > s} \sup_{j \leq J_{t-1}} |e^{2\pi ia_jm_t} - 1| \\
&\leq \sum_{t > s} \frac{1}{50r} < \frac{1}{50}.
\end{split}
\end{equation*}

\noindent Hence, we have that for all $s$

\begin{equation*}
\left| J_s^{-1} \sum_{j \leq J_s} e^{2\pi ia_j n(\bar{\alpha})} -
J_s^{-1} \sum_{j \leq J_s} e^{2\pi i a_jn_s} \right| < \frac{1}{40}
+ \frac{1}{50} < \frac{1}{20}.
\end{equation*}

\noindent Now, if $\alpha_s = 0$, then $n_s = 0$ and

\begin{equation*}
\left| J_s^{-1} \sum_{j \leq J_s} e^{\pi ia_j n(\bar{\alpha})} - 1
\right| < 1/20 < 1/10.
\end{equation*}

\noindent If $\alpha_s = 1$, then $n_s = m_s$ and

\begin{equation*}
\begin{split}
\left| J_s^{-1} \sum_{j \leq J_s} e^{\pi ia_j n(\bar{\alpha})} - 1
\right| &\geq \left| J_s^{-1} \sum_{j \leq J_s} e^{\pi ia_j m_s} - 1
\right| - \left| J_s^{-1} \sum_{j \leq J_s} e^{\pi ia_j
n(\bar{\alpha})} - J_s^{-1} \sum_{j
\leq J_s} e^{\pi ia_j m_s} \right| \\
&> 3/4 - 1/20 > 1/2.
\end{split}
\end{equation*}

\noindent This completes the proof.
\end{proof}

\vskip0.2in

\begin{proof}[Proof of Lemma~\ref{lemma:comb}]
Suppose not.  Then, there exist subsets $M_1, M_2 \subset \{1,
\ldots, r\}$ with $|M_1|, |M_2| < r_0$ so that for any $j \notin
M_1$ there is an $m_1 \in M_1$ such that $\|\re h_j - \re h_{m_1}\|
\leq \alpha/4$, and for any $j \notin M_2$ there is an $m_2 \in M_2$
such that $\|\im  h_j - \im  h_{m_2}\| \leq \alpha/4$. Define $M =
M_1 \cup M_2$.  Then, $|M| = m < 2r_0$.

For each $j \notin M$, we can associate a pair $(a,b)$ with $a, b
\in M$ such that $\|\re  h_j - \re  h_a\| \leq \alpha/4$ and $\|\im
\, h_j - \im  h_b\| \leq \alpha/4$.  Associate one such pair to each
$j \notin M$. Now, there are  $m^2$ distinct such pairs. But, there
are $r - m$ points $j \notin M$.  As $m < 2r_0$ and by the defintion
of $r$, we have $r - m > m^2$.  Therefore, two distinct points $j, k
\notin M$ must be assigned the same pair, say $(a,b)$. Namely,

\begin{equation*}
\begin{split}
\|h_j - h_k\| &\leq \|\re  h_j - \re  h_k\| + \|\im  h_j - \im  h_k\| \\
&\leq \|\re  h_j - \re  h_a\| + \|\re  h_a - \re  h_k\| + \|\im h_j
- \im h_b\| + \|\im  h_b - \im  h_k\| \\
&\leq \alpha.
\end{split}
\end{equation*}

\noindent This contradicts the hypothesis.
\end{proof}

\section{Comments}

\begin{enumerate}

\item The only significant difference between the proof of the first
entropy result here and in~\cite{bourgain} is the use of
Theorem~\ref{thm:fernique}.  Bourgain uses a different statement,
namely that there is some $c > 0$ so that

\begin{equation*}
P \left\{ \omega : \mu \left\{x : F^*(x,\omega) \geq c \int_\Omega
F^*(x,\omega') \, P(d\omega') \right\} > c \right\} > c
\end{equation*}

\noindent for all $N, J \in \mathbb{N}$ big enough and $f \in
L^\infty(\mu)$, $\|f\|_{L^2(\mu)} \leq 1$.
Theorem~\ref{thm:fernique} is used in an almost identical way to
this statement in the proofs of both entropy results.

\item Aside from the above comment, the proofs of the second entropy
result here and in~\cite{bourgain} differ in only one other place,
and only slightly.  Bourgain states and uses a different Banach
principle for $L^\infty(\mu)$.  In particular, he states that if
$S_n f$ converges almost surely for all $f \in L^\infty(\mu)$, then
there is some function $\delta(\epsilon)$, which goes to 0 as
$\epsilon \rightarrow 0$, such that $\int_X \sup_n |S_n f| \, d\mu <
\delta(\epsilon)$ whenever $\|f\|_{L^\infty(\mu)} \leq 1$ and
$\|f\|_{L^1(\mu)} < \epsilon$.  This result can be used in the same
manner as Theorem~\ref{thm:bj2}.  I have included the result from
Bellow and Jones instead, because their proof is so easily
understood.

\item For more on this and related topics, see work by Roger
Jones~\cite{jones1, jones2}, Michael Lacey~\cite{lacey}, and Michel
Weber (with Mikhail Lifshits and Dominique Schneider) [8 - 19].

\item I will be happy to field any questions or concerns via e-mail.  I
would also appreciate being alerted to any typos.  Be sure to use a
descriptive subject, as I get lots of spam.

\end{enumerate}

\end{document}